\documentclass[11pt]{article}

\usepackage{amsmath, amssymb, bm, color}
\usepackage{times,mathptmx}

\newtheorem{Thm}{Theorem}
\newtheorem{Lem}{Lemma}
\newtheorem{Prop}{Proposition}

\newtheorem{Rem}{Remark}
\newtheorem{Exa}{Example}

\begin{document}
{\allowdisplaybreaks{

\title{
Joint distributions for stochastic functional differential equations\footnotetext[0]{
To appear in {\it{Stochastics: An International Journal of Probability and Stochastic Processes}}}\footnote{Dedicated to Professor Takashi Komatsu on the occasion of his 70th birthday. }
}
\author{
{\sc{Atsushi Takeuchi}}\footnote{
E-mail address: {\texttt{takeuchi@sci.osaka-cu.ac.jp}}. \ 
Postal address: Department of Mathematics, Osaka City University, 
Sugimoto 3-3-138, Sumiyoshi-ku, Osaka 558-8585, Japan.
}
}
\date{}
\maketitle

\begin{abstract}
Consider stochastic functional differential equations, whose coefficients depend on past histories. 
The solution determines a non-Markov process. 
In the present paper, we shall obtain the existence of smooth densities for joint distributions of solutions, under the uniformly elliptic condition on the diffusion coefficients, via the Malliavin calculus. 
As an application, we shall study the computations of the Greeks on options associated with the asset price dynamics models with delayed effects. 

\vspace{5pt}
\end{abstract}

\section{Introduction}\label{Section 1}
Let $(\Omega, \mathcal{F}, \mathbb{P})$ be a probability space, and $r$ and $T$ be positive constants, which are fixed throughout the paper. 
Denote by $W = \big\{ W(t) = \big( W^1(t), \dots, W^m(t) \big) \, ; \, 0 \le t \le T \big\}$ the $m$-dimensional Brownian motion starting from the origin. 
Write $\mathcal{F}_t = \sigma \big[ W(s) \, ; \, 0 \le s \le t \big] \vee \mathcal{N}$ for $0 \le t < T$ and $\mathcal{F}_T = \mathcal{F}$, where $\mathcal{N}$ is the family of $\mathbb{P}$-null sets. 
Let $A_0, \, A_1, \, \dots, \, A_m$ be in $C_{1+,b}^{1, \infty} \big( [0,T] \times C([-r,0] \, ; \, \mathbb{R}^d) \, ; \, \mathbb{R}^d)$, that is, those are jointly continuous in $(t,f) \in [0,T] \times C([-r,0] \, ; \, \mathbb{R}^d)$ such that, for each $0 \le i \le m$, 
\begin{itemize}
\item for each $f \in C([-r,0] \, ; \, \mathbb{R}^d)$, the mapping $[0,T] \ni t \longmapsto A_i (t, f) \in \mathbb{R}^d$ is differentiable such that its partial derivative $\partial_t \, A_i (t,f)$ is bounded, 
\item for each $t \in [0,T]$, the mapping $C([-r,0] \, ; \, \mathbb{R}^d) \ni f \longmapsto A_i (t,f) \in \mathbb{R}^d$ is smooth in the Fr\'echet sense such that all partial Fr\'echet derivatives $\nabla^k \, A_i (t,f) \ (k \in \mathbb{N})$ are bounded. 
\end{itemize}
For a deterministic path $\eta \in C([-r,0] \, ; \, \mathbb{R}^d)$, consider {\it{the stochastic functional differential equation}}: 
\begin{equation}\label{eq:SFDE X}
X(t) = 
\begin{cases}
\eta(t) 
& (-r \le t \le 0), \\
\displaystyle \eta(0) + \int_0^t A_0 (s, X_s) \, \mbox{d}s 
+ \int_0^t A (s, X_s) \, \mbox{d}W(s) 
& (0 \le t \le T), 
\end{cases}
\end{equation}
where $A = \big( A_1, \dots, A_m \big)$, and $X_t = \big\{ X(t+u) \, ; \, -r \le u \le 0 \big\}$ is the segment. 
Since the coefficients $A_0, \, A_1, \, \dots, \, A_m$ satisfy the Lipschitz condition and the linear growth one, there exists a unique solution to \eqref{eq:SFDE X} such that 
$$
\mathbb{E} \left[ \sup_{-r \le t \le T} |X(t)|^p \right] \le C_{1, p, \eta, T}
$$
for any $p > 1$. 
Moreover, the solution process $X  = \big\{ X(t) \, ; \, -r \le t \le T \big\}$ is non-Markovian, because the current state of the process depends on the whole past histories of the process $X$. 
See \cite{IN, Mohammed1, Mohammed2}. 
Thus, we cannot use any methods in analysis, partial differential equations and potential theory at all. 

Such equation was initiated by It\^o and Nisio \cite{IN} more than 50 years ago. 
It seems us to be very natural to study the models described by \eqref{eq:SFDE X} stated above, because the models with their past histories often appear in finance, physics, biology and industry, etc. 
One of the typical examples in mathematical finance is {\it{the delayed Black-Scholes model}} studied in \cite{Arriojas-Hu-Mohammed-Pap, Chang-Pang-Pemy, Chang-Youree-1, Chang-Youree-2, MS}, which will be mentioned in Section \ref{Section 5}. 
On the other hand, the Malliavin calculus is applicable to the study on the densities for the solution to \eqref{eq:SFDE X}. 
Kusuoka and Stroock \cite{KS} obtained the result on the existence of the smooth density for the solution with respect to the Lebesgue measure under the uniformly elliptic condition on the diffusion coefficients $A_1, \, \dots, \, A_m$. 
Bell and Mohammed in \cite{Bell-Mohammed-1, Bell-Mohammed-2} also studied the same problem in case of {\it{stochastic delay differential equations }}such that $A_i (t,f) = \hat{A}_i (t,f(-r)) \ (i = 1 \, \dots, m)$ for $t \in [0,T]$ and $f \in C([-r,0] \, : \mathbb{R}^d)$, where $\hat{A}_i: [0,T] \times \mathbb{R}^d \to \mathbb{R}^d$ with some conditions on the boundedness and the regularity. 
They obtained in \cite{Bell-Mohammed-1, Bell-Mohammed-2} the existence of the smooth density under the degeneracy condition on $\sum_{i=1}^m \hat{A}_i \, \hat{A}_i^\ast$ by using the delay structure of the equation and conditioning the past history of the process, which is weaker than the uniformly elliptic condition on $\sum_{i=1}^m A_i \, A_i^\ast$. 
Furthermore, Kitagawa and Takeuchi \cite{KT} studied the asymptotic behavior of the density such as the Varadhan-type estimate for diffusion processes, by the large deviation theory and the Malliavin calculus, in which the constant $r$, called {\it{the delay parameter}}, plays a crucial role. 

In the present paper, we will study the finite-dimensional joint distribution on the solution process to \eqref{eq:SFDE X}, from the viewpoint of the Malliavin calculus. 
As stated in Theorem \ref{Thm:n-state} below, the joint distribution of the solution admits a smooth density under the uniformly elliptic condition on the diffusion coefficients $A_1, \, \dots, \, A_m$. 
As an application, we shall also study the sensitivity analysis on the solution with respect to the initial state, which can be regarded as the computation of the Greeks for the options on the delayed asset price dynamics model. 

The paper is organized as follows: 
Section \ref{Section 2} is devoted to a brief introduction of the Malliavin calculus and its application to stochastic functional differential equations. 
The result on the existence of the smooth density for the finite-dimensional joint distribution associated with the solution will be stated in Section \ref{Section 3}. 
Sections \ref{Section 4} and \ref{Section 5} are typical applications of our result. 
In Section \ref{Section 4}, we will study the sensitivity of the discrete and integral average related to the solution. 
The key points are to give the estimates on the Malliavin covariance matrices, which are implied by the uniformly elliptic conditions on the diffusion coefficients $A_1, \, \dots, \, A_m$ of \eqref{eq:SFDE X}. 
We will study a delayed Black-Scholes model raised in \cite{Arriojas-Hu-Mohammed-Pap, MS} in Section \ref{Section 5}, in order to compute the Greeks on the options. 

\section{Malliavin calculus}\label{Section 2}

In this section, we shall apply the Malliavin calculus to the stochastic functional differential equation \eqref{eq:SFDE X}. 
See \cite{Nualart} on details of the Malliavin calculus. 
Let $\big( \mathbb{W}_0^m, \mathcal{W}, \mathbb{P}^{\mathbb{W}_0^m} \big)$ be the Wiener space, that is, $\mathbb{W}_0^m$ is the set of $\mathbb{R}^m$-valued continuous functions on $[0,T]$ starting from the origin in $\mathbb{R}^m$, $\mathcal{W}$ is the topological $\sigma$-algebra on $\mathbb{W}_0^m$, and $\mathbb{P}^W$ is the Wiener measure over $(\mathbb{W}_0^m, \mathcal{W})$. 
Denote by $\mathbb{H}_0^m$ be the Cameron-Martin subspace of $\mathbb{W}_0^m$ with the inner product 
$$
\langle g, h \rangle_{\mathbb{H}_0^m}
:= \int_0^T \langle \dot{g} (u), \dot{h} (u) \rangle_{\mathbb{R}^m} \, \mbox{d}u \quad ( g, \, h \in \mathbb{H}_0^m), 
$$
where $\dot{g}(u)$ is the derivative of $g$ in $u$. 

For $0 \le s \le T$, let $ \big\{ Z(t, s) \, ; \, -r \le t \le T \big\}$ be the $\mathbb{R}^d \otimes \mathbb{R}^d$-valued process determined by
\begin{equation}\label{eq:Z}
Z(t,s) = 
\begin{cases}
\bm{0} \qquad \qquad \qquad \qquad \qquad \qquad \qquad \qquad \qquad \qquad \qquad (-r \le t \le 0 \ \mbox{or} \ t < s), \\
\displaystyle I_d + \int_s^t \nabla A_0 (u, X_u) \, Z_u (\cdot , s) \, \mbox{d}u 
+ \int_s^t \sum_{i=1}^m \nabla A_i (u, X_u) \, Z_u (\cdot, s) \, \mbox{d}W^i (u) \quad (s \le t \le T), 
\end{cases}
\end{equation}
where $I_d \in \mathbb{R}^d \otimes \mathbb{R}^d$ is the identity, and $Z_u (\cdot, s) = \big\{ Z(u + \tau, s) \, ; \, -r \le \tau \le 0 \big\}$. 
Then, we have 
\begin{Prop}[cf. \cite{KT, KS}]\label{Prop:Malliavin smoothness of X(t)}
For $-r \le t \le T$, the random variable $X(t)$ is smooth in the Malliavin sense. 
Moreover, the Malliavin derivative $DX(t) = \big\{ D_u X(t) \, ; \, 0 \le u \le T \big\}$ of $X(t)$ and the Malliavin covariance matrix $V(t) := \langle D X(t), DX(t) \rangle_{\mathbb{H}_0^m}$ for $X(t)$ can be computed as follows: 
\begin{gather}
D_u X(t) = \int_0^{u \wedge t} Z(t, s) \, A(s, X_s) \, \mbox{d}s, \label{eq:computation of DX(t)} \\
V(t) = \int_0^t \sum_{i=1}^m Z(t, s) \, A_i(s, X_s) \, \big\{ Z(t, s) \, A_i(s, X_s) \big\}^\ast \, \mbox{d}s. \label{eq:MC for X(t)}
\end{gather}
\end{Prop}
{\it{Proof}}. 
The Malliavin smoothness of $X(t)$ can be justified by the limiting argument via the successive approximation $X^{(n)} = \big\{ X^{(n)}(t) \, ; \, -r \le t \le T \big\} \ (n \in \mathbb{N} \cup \{ 0 \})$ of the process $X$: 
\begin{align*}
X^{(0)} (t) & = \eta (t) \, \mathbb{I}_{[-r,0]}(t) + \eta(0) \, \mathbb{I}_{(0,T]} (t), \\
X^{(n+1)}(t) & = 
\begin{cases}
\eta (t) & (-r \le t \le 0), \\
\displaystyle \eta (0) + \int_0^t A_0 (s, X_s^{(n)}) \, \mbox{d}s 
+ \int_0^t A (s, X_s^{(n)}) \, \mbox{d}W(s) & (0 < t \le T) 
\end{cases}
\end{align*}
for $n \in \mathbb{N} \cup \{ 0 \}$, and the inductive argument on the order of the derivatives. 
On the other hand, since $D_u X^{(n)} (t) = \bm{0}$ for $-r \le t \le 0$ and $n \in \mathbb{N} \cup \{ 0 \}$, and 
\begin{align*}
D_u X^{(0)} (t) 
& = \bm{0}, \\
D_u X^{(n)}(t) 
& = \int_0^{u \wedge t} A(s, X_s^{(n-1)}) \, \mbox{d}s 
+ \int_0^t \nabla A_0 (s, X_s^{(n-1)}) \, D_u X_s^{(n-1)} \mbox{d}s \\
& \qquad + \int_0^t \sum_{i=1}^m \nabla A_i (s, X_s^{(n-1)}) \, D_u X_s^{(n-1)} \mbox{d}W^i(s) 
\end{align*}
for $n \in \mathbb{N}$ and $0 \le t \le T$, the limiting argument leads us to see that 
\begin{equation}\label{eq:DX(t)}
D_u X(t) = 
\begin{cases}
\bm{0} \qquad \qquad \qquad \qquad \qquad \qquad \qquad \qquad \qquad \qquad \qquad \qquad \qquad \quad (-r \le t \le 0), \\[5pt]
\displaystyle \int_0^{u \wedge t} \!\! A(s, X_s) \mbox{d}s 
+ \int_0^t \nabla A_0 (s, X_s) D_u X_s \mbox{d}s 
+ \int_0^t \sum_{i=1}^m \nabla A_i (s, X_s) D_u X_s \mbox{d}W^i(s) \  (u \le t \le T). 
\end{cases}
\end{equation}
Thus, we can derive \eqref{eq:computation of DX(t)}, because of the uniqueness of the solution to the equation \eqref{eq:DX(t)}. 
Moreover, the Malliavin covariance matrix $V(t)$ can be computed directly as \eqref{eq:MC for X(t)}. 
\hfill $\square$

\section{Density of joint distributions}\label{Section 3}

In this section, we shall mention the main result of the present paper. 
Before doing that, recall the classical result on the existence of the smooth density for the probability law of $X(t)$. 

\begin{Lem}[cf. \cite{KT, KS}]\label{Lem:density for X(t)}
Suppose that the coefficients $A_1, \, \dots, \, A_m$ satisfy the uniformly elliptic condition: 
there exists a constant $C_2 > 0$ such that 
\begin{equation}\label{UE condition}
\inf_{0 \le t \le T} \, \inf_{f \in C([-r,0] \, ; \, \mathbb{R}^d)} \ \inf_{v \in \mathbb{S}^{d-1}} \sum_{i=1}^m \big( A_i (t, f) \cdot v \big)^2 
\ge C_2. 
\end{equation}
Then, for each $0 < t \le T$, the probability law of the $\mathbb{R}^d$-valued random variable $X(t)$ admits a smooth density with respect to the Lebesgue measure on $\mathbb{R}^d$. 
\end{Lem}

Let $n \in \mathbb{N}$ be arbitrary, and $0 < t_1 < \dots < t_{n-1} < t_n = t$. 
Before introducing our result, we shall give an easy example, which is our motivation of our interests in the present paper. 
\begin{Exa}
{\rm{
Consider the case of $m = d$, $A_0 (t,f) \equiv 0$, and $A (t,f) = I_d$. 
Let $t_0 = 0$, and write 
$$
p(t, x, y) = \frac{1}{\sqrt{2 \pi t}} \, \exp \left( - \frac{ (y-x)^2}{2t} \right). 
$$
Then, the process $X$ is the $d$-dimensional Brownian motion. 
Since 
\begin{align*}
\mathbb{P} \big[ W(t_1) \in K_1, \, \dots, \, W(t_n) \in K_n \big] 
= \int_{K_1 \times \cdots \times K_n} p(t_1, 0, y_1) \, \prod_{k=2}^n p(t_k -t_{k-1}, y_{k-1}, y_k) \, \mbox{d}y_1 \cdots \mbox{d}y_n
\end{align*}
for $K_1, \, \dots, \, K_n \in \mathcal{B}(\mathbb{R}^d)$, the joint distribution $\big( W(t_1), \dots, W(t_n) \big)$ is absolutely continuous with respect to the Lebesgue measure over $\mathbb{R}^{nd}$ such that its density function 
$$
\varphi_{t_1, \dots, t_n} (y_1, \dots, y_n) 
:= p(t_1, 0, y_1) \, \prod_{k=2}^n p(t_k -t_{k-1}, y_{k-1}, y_k)
$$
is smooth in $(y_1, \dots, y_n) \in \mathbb{R}^{nd}$. 
\hfill $\square$
}}
\end{Exa}

Now, we shall introduce the result in this paper, which can be regarded as the natural extension of Lemma \ref{Lem:density for X(t)}. 

\begin{Thm}\label{Thm:n-state}
Under the condition \eqref{UE condition} in Lemma \ref{Lem:density for X(t)}, the joint distribution of $\big( X(t_1), \dots, X(t_n) \big)$ admits a smooth density with respect to the Lebesgue masure on $\mathbb{R}^{nd}$. 
\end{Thm}
\noindent
{\it{Proof}}. 
The $\mathbb{R}^{nd}$-valued random variable $\big( X(t_1), \dots, X(t_n) \big)$ determined by the equation \eqref{eq:SFDE X} is smooth in the Malliavin sense, as stated in Proposition \ref{Prop:Malliavin smoothness of X(t)}, because so are all of $X(t_k) \ (1 \le k \le n)$. 
Moreover, the corresponding Malliavin covariance matrix $V(t_1, \dots, t_n)$ for $\big( X(t_1), \dots, X(t_n) \big)$ is 
\begin{align*}
V(t_1, \dots, t_n) 
& = 
\begin{pmatrix}
\langle DX(t_1), DX(t_1) \rangle_{\mathbb{H}_0^m} & \cdots & \langle DX(t_1), DX(t_n) \rangle_{\mathbb{H}_0^m} \\
\vdots & \ddots & \vdots \\
\langle DX(t_n), DX(t_1) \rangle_{\mathbb{H}_0^m} & \cdots & \langle DX(t_n), DX(t_n) \rangle_{\mathbb{H}_0^m}
\end{pmatrix}
\\[5pt]
& = 
\begin{pmatrix}
\displaystyle \int_0^T \Phi(t_1,s) \, \Phi(t_1,s)^\ast \, \mbox{d}s 
& \cdots & 
\displaystyle \int_0^T \Phi(t_1,s) \, \Phi(t_n,s)^\ast \, \mbox{d}s \\[10pt]
\vdots & \ddots & \vdots \\[10pt]
\displaystyle \int_0^T \Phi(t_n,s)  \, \Phi(t_1,s)^\ast \, \mbox{d}s 
& \cdots & 
\displaystyle \int_0^T \Phi(t_n,s) \, \Phi(t_n,s)^\ast \, \mbox{d}s
\end{pmatrix}
\end{align*}
where $\Phi(t,s) := Z(t,s) \, A(s, X_s) \, \mathbb{I}_{(s \le t)}$ is $\mathbb{R}^{nd} \otimes \mathbb{R}^{nd}$-valued. 

All we have to do is to study the negative-order moment of $\det V(t_1, \dots, t_n)$. 
Let $v = \big( v_1, \dots, v_n \big) \in \mathbb{R}^{nd}$ such that $|v| = \big(|v_1|^2 + \cdots + |v_n|^2 \big)^{1/2} = 1$, and write $t_0 = 0$. 
Then, we have 
\begin{align*}
\left\langle v, V(t_1, \dots, t_n) \, v \right\rangle_{\mathbb{R}^{nd}} 
& = \int_0^{t_n} \left| \sum_{k=1}^n \Phi(t_k,s)^\ast v_k \right|^2 \mbox{d}s \\
& = \sum_{j=1}^n \int_{t_{j-1}}^{t_j} \left| \sum_{k=j}^n \Phi(t_k,s)^\ast v_k \right|^2 \mbox{d}s \\
& \ge \sum_{j=1}^{n-1} \int_{t_{j-1}}^{t_j} \left| \sum_{k=j}^n \Phi(t_k,s)^\ast v_k \right|^2 \mbox{d}s \, \mathbb{I}_{(v_n = 0)}
+ \int_{t_{n-1}}^{t_n} \left| \Phi(t_n,s)^\ast v_n \right|^2 \mbox{d}s \, \mathbb{I}_{(v_n \neq 0)} \\
& \ge  \sum_{j=1}^{n-2} \int_{t_{j-1}}^{t_j} \left| \sum_{k=j}^n \Phi(t_k,s)^\ast v_k \right|^2 \mbox{d}s \, \mathbb{I}_{(v_{n-1} = 0, v_n = 0)} \\
& \quad + \int_{t_{n-2}}^{t_{n-1}} \left| \Phi(t_{n-1},s)^\ast v_{n-1} \right|^2 \mbox{d}s \, \mathbb{I}_{(v_{n-1} \neq 0, v_n = 0)}
+ \int_{t_{n-1}}^{t_n} \left| \Phi(t_n,s)^\ast v_n \right|^2 \mbox{d}s \, \mathbb{I}_{(v_n \neq 0)} \\
& \ge \cdots \\
& \ge \sum_{j=1}^{n-1} \int_{t_{j-1}}^{t_j} \left| \Phi(t_j,s)^\ast v_j \right|^2 \mbox{d}s \, \mathbb{I}_{(v_j \neq 0, v_{j+1} = 0, \dots, v_n = 0)} 
+ \int_{t_{n-1}}^{t_n} \left| \Phi(t_n,s)^\ast v_n \right|^2 \mbox{d}s \, \mathbb{I}_{(v_n \neq 0)}.  
\end{align*}
Remark that 
\begin{align*}
I_j 
& := \int_{t_{j-1}}^{t_j} \left| \Phi(t_j,s)^\ast v_j \right|^2 \mbox{d}s \, \mathbb{I}_{(v_j \neq 0, v_{j+1} = 0, \dots, v_n = 0)} \\
& = \int_{t_{j-1}}^{t_j} \sum_{i=1}^m \left\langle Z(t_j, s) \, A_i(s, X_s), v_j \right\rangle_{\mathbb{R}^d}^2 \, \mbox{d}s \ \mathbb{I}_{(v_j \neq 0, v_{j+1} = 0, \dots, v_n = 0)} \\
& \ge \int_{t_j - \lambda^{-\alpha}}^{t_j} \sum_{i=1}^m \left\langle Z(t_j, s) \, A_i(s, X_s), v_j \right\rangle_{\mathbb{R}^d}^2 \, \mbox{d}s \ \mathbb{I}_{(v_j \neq 0, v_{j+1} = 0, \dots, v_n = 0)} \\
& \ge \frac{ \ 1 \ }{2} \int_{t_j - \lambda^{-\alpha}}^{t_j} \sum_{i=1}^m \left\langle A_i(s, X_s), v_j \right\rangle_{\mathbb{R}^d}^2 \, \mbox{d}s \ \mathbb{I}_{(v_j \neq 0, v_{j+1} = 0, \dots, v_n = 0)} 
- \int_{t_j - \lambda^{-\alpha}}^{t_j} \sum_{i=1}^m | A_i(s, X_s) |^2 |v_j|^2 \| Z(t_j,s) - I_d \|^2 \mbox{d}s \\
& \ge \frac{C_2 \, \lambda^{-\alpha}}{2} |v_j|^2 \, \mathbb{I}_{(v_j \neq 0, v_{j+1} = 0, \dots, v_n = 0)} - C_3^2 \, \lambda^{-\beta}
\end{align*}
on the subset 
$$
\Omega_2 
:= \bigcap_{j=1}^n \left( \left\{ \sup_{t_j - \lambda^{-\alpha} \le s \le t_j} |X(s)| \le C_3 \right\} 
\cap \left\{ \int_{t_j - \lambda^{-\alpha}}^{t_j} \big\| Z(t_j,s) - I_d \big\|^2 \mbox{d}s \le \lambda^{- \beta} \right\} \right)
$$
under the condition \eqref{UE condition} on $A_1, \, \dots, \, A_m$, where $0 < \alpha < 1$ and $\alpha < \beta < 2 \alpha$. 
Hence, it holds that 
\begin{align*}
\left\langle v, V(t_1, \dots, t_n) \, v \right\rangle_{\mathbb{R}^{nd}} 
& \ge \frac{C_2 \, \lambda^{-\alpha}}{2} \left( \sum_{j=1}^{n-1} |v_j|^2 \, \mathbb{I}_{(v_j \neq 0, v_{j+1} = 0, \dots, v_n = 0)} + |v_j|^2 \, \mathbb{I}_{(v_n \neq 0)} \right)
- n \, C_3^2 \, \lambda^{-\beta}
\end{align*}
on $\Omega_2$. 
On the other hand, we shall remark that 
$$
\mathbb{P} \big[ \Omega_2^c \big] 
\le n \left( C_3^{-p} \, C_{1,p,\eta,T} \, \lambda^{-\alpha p /2} + C_{4,p,T} \, \lambda^{- (2 \alpha - \beta)p} \right)
$$
for any $p > 1$ from the Chebyshev inequality. 

Now, we shall return to study the upper estimate of 
$$
I(\lambda) := \sup_{ v = (v_1, \dots, v_n) \in \mathbb{R}^{nd}, |v| = 1} \mathbb{E} \left[ \exp \big( - \lambda \, \left\langle v, V(t_1, \dots, t_n) \, v \right\rangle_{\mathbb{R}^{nd}} \big) \right]. 
$$
Since the mapping $\mathbb{S}^{nd-1} \ni v \longmapsto \mathbb{E} \left[ \exp \big( - \lambda \, \left\langle v, V(t_1, \dots, t_n) \, v \right\rangle_{\mathbb{R}^{nd}} \big) \right] \in \mathbb{R}$ is continuous, we can find 
$$
\tilde{v} = \mbox{argmax} \left\{ \mathbb{E} \left[ \exp \big( - \lambda \, \left\langle v, V(t_1, \dots, t_n) \, v \right\rangle_{\mathbb{R}^{nd}} \big) \right] \, ; \, v = (v_1, \dots, v_n) \in \mathbb{R}^{nd}, |v| = 1 \right\}. 
$$
Therefore, we can conclude that 
\begin{align*}
I(\lambda) 
& \equiv \mathbb{E} \left[ \exp \big( - \lambda \, \left\langle \tilde{v}, V(t_1, \dots, t_n) \, \tilde{v} \right\rangle_{\mathbb{R}^{nd}} \big) \right] \\
& \le \exp \left[ - \lambda \left(\frac{C_2 \, \lambda^{-\alpha}}{2} \left( \sum_{j=1}^{n-1} |\tilde{v}_j|^2 \, \mathbb{I}_{(\tilde{v}_j \neq 0, \tilde{v}_{j+1} = 0, \dots, \tilde{v}_n = 0)} + |\tilde{v}_j|^2 \, \mathbb{I}_{(\tilde{v}_n \neq 0)} \right)
- n \, C_3^2 \, \lambda^{-\beta} \right) \right] \\
& \qquad + n \left( C_3^{-p} \, C_{1,p,\eta,T} \, \lambda^{- \alpha p /2} + C_{4,p,T} \, \lambda^{- (2 \alpha - \beta) p} \right) \\
& = o(\lambda^{- C_5 \, p})
\end{align*}
as $\lambda \to +\infty$ for any $p > 1$, so we can obtain 
\begin{align*}
\mathbb{E} \big[ \big( \det V(t_1, \dots, t_n) \big)^{-q} \big] 
& = \mathbb{E} \left[ \Big( \inf_{v \in \mathbb{S}^{nd-1}} \left\langle v, V(t_1, \dots, t_n) \, v \right\rangle_{\mathbb{R}^{nd}} \Big)^{-nqd} \right] \\
& \le C_6 \, \sup_{v \in \mathbb{S}^{nd-1}} \mathbb{E} \Big[ \big( \left\langle v, V(t_1, \dots, t_n) \, v \right\rangle_{\mathbb{R}^{nd}} \big)^{- (nqd + 4nd -4) } \Big] + C_7 \\
& = \frac{C_6}{\Gamma (nqd + 4nd -4)} \int_0^{+\infty} \lambda^{4qd + 4nd -5} \, I(\lambda) \, \mbox{d}\lambda +C_7 
< +\infty 
\end{align*}
for any $q > 1$. 
See \cite{KomatsuTakeuchi} on the detailed discussion of the second inequality stated above. 
\hfill $\square$

\begin{Rem}
{\rm{
Bell and Mohammed \cite{Bell-Mohammed-2} have studied {\it{the stochastic delay equation with hereditary drift}}: 
\begin{equation}\label{Eq:B-M}
X(t) = 
\begin{cases}
\eta(t) & (-r \le t \le 0), \\
\displaystyle \eta(0) + \int_0^t \hat{A}_0 \big( s, \{ X(u) \, ; \, -r \le u \le s \} \big) \, \mbox{d}s + \int_0^t \hat{A} \big( s, X(s-r) \big) \, \mbox{d}W(s) & (0 \le t \le T), 
\end{cases}
\end{equation}
where $\hat{A}_0 : [0,T] \times C \big( [-r,T] \, ; \, \mathbb{R}^d \big)$ such that $\hat{A}_0 (t,f) \ \big( t \in [0,T], \, f \in C ( [-r,T] \, ; \, \mathbb{R}^d ) \big)$ depends only on $\big\{ f(s) \, ; \, -r \le s \le t \big\}$, and $\hat{A}_i : [0,T] \times \mathbb{R}^d \to \mathbb{R}^d \ (i = 1, \dots, m)$ with the certain conditions on the boundedness and the regularity. 
Write $\hat{A} = \big( \hat{A}_1, \dots, \hat{A}_m \big)$, and let $t \in [0,T]$. 
They showed in \cite{Bell-Mohammed-2} the existence of the smooth density for the law of $X(t)$, under the degeneracy of the $\mathbb{R}^d \otimes \mathbb{R}^d$-valued function $\hat{A} \, \hat{A}^\ast $ of polynomial order on hypersurfaces in $\mathbb{R}^d$, which is weaker than the condition \eqref{UE condition} in the present paper, by using the delay structure in \eqref{Eq:B-M} and conditioning on the past history of the process. 

As for the equation \eqref{Eq:B-M}, we can also derive the same assertions as Lemma \ref{Lem:density for X(t)} and Theorem \ref{Thm:n-state}, under the degeneracy condition on $\hat{A} \, \hat{A}^{\ast}$ as stated in  \cite{Bell-Mohammed-2}. 
In fact, since it can be checked, similarly to Proposition \ref{Prop:Malliavin smoothness of X(t)},  that the Malliavin calculus is applicable to the solution of the equation \eqref{Eq:B-M}, our goal in the argument is to study the negative-order moment of the determinants on the corresponding Malliavin covariance matrices for $X(t)$ and $\big( X(t_1), \dots, X(t_n) \big)$ of the forms: 
\begin{align}
V(t) 
& = \int_{0}^{t} D_{u}X(t) \, \big( D_{u} X(t) \big)^{\ast} \, \mbox{d}u, \\[5pt]
V(t_1, \dots, t_{n}) 
& = 
\begin{pmatrix}
\displaystyle \int_{0}^{t_{n}} D_{u}X(t_{1}) \, \big( D_{u} X(t_{1}) \big)^{\ast} \, \mbox{d}u 
& \cdots 
& \displaystyle \int_{0}^{t_{n}} D_{u}X(t_{1}) \, \big( D_{u} X(t_{n}) \big)^{\ast} \, \mbox{d}u \\
\vdots 
& \ddots 
& \vdots \\
\displaystyle \int_{0}^{t_{n}} D_{u}X(t_{n}) \, \big( D_{u} X(t_{1}) \big)^{\ast} \, \mbox{d}u 
& \cdots 
& \displaystyle \int_{0}^{t_{n}} D_{u}X(t_{n}) \, \big( D_{u} X(t_{n}) \big)^{\ast} \, \mbox{d}u
\end{pmatrix}
,
\end{align}
where $0 = t_{0} < t_{1} < \cdots < t_{n-1} < t_{n} = t$ such that $\| \Delta \| := \max_{1 \le j \le n} ( t_{j} - t_{j-1} ) \le r$. 
Then, we have only to check the lower estimates of the quadratic forms of $V(t)$ and $V(t_{1}, \dots, t_{n})$: 
\begin{align}
\langle v, V(t) \, v \rangle_{\mathbb{R}^{d}}
& = \int_{0}^{t} \big| (D_{u} X(t) \big)^{\ast} \, v \big|^{2} \, \mbox{d}u, \\
\langle (v_{1}, \dots, v_{n}), V(t_{1}, \dots, t_{n}) \, (v_{1}, \dots, v_{n}) \rangle_{\mathbb{R}^{{nd}}}
& \ge \sum_{j=1}^{n-1} \int_{t_{j-1}}^{t_{j}} \big| (D_{u} X(t_{j}))^{\ast} \, v_{j} \big|^{2} \, \mbox{d}u \, \mathbb{I}_{(v_{j} \neq 0, v_{j+1} = 0, \dots, v_{n}= 0)} \nonumber \\
& \qquad + \int_{t_{n-1}}^{t_{n}} \big| (D_{u} X(t_{n}))^{\ast} \, v_{n} \big|^{2} \, \mbox{d}u \, \mathbb{I}_{(v_{n} \neq 0)}, 
\end{align}
where $v \in \mathbb{R}^{d}$ with $|v| = 1$, and $(v_{1}, \dots, v_{n}) \in \mathbb{R}^{nd}$ with $\sum_{j=1}^{n} |v_{j}|^{2} = 1$. 

Finally, we shall remark that it would be open whether Lemma \ref{Lem:density for X(t)} and Theorem \ref{Thm:n-state} on the equation \eqref{eq:SFDE X} can be obtained, under the degeneracy condition as stated in \cite{Bell-Mohammed-2}, because the special forms of the diffusion coefficients $\hat{A}_{1}, \, \dots, \, \hat{A}_{m}$ in the equation \eqref{Eq:B-M} play a crucial role in \cite{Bell-Mohammed-2},  
\hfill $\square$
}}
\end{Rem}

\section{Applications}\label{Section 4}

In this section, we shall study one of the typical applications of Theorem \ref{Thm:n-state}. 
Consider the case of $d = m = 1$ throughout this section, in order to avoid the complication of our argument. 
Let $f \in C^1( \mathbb{R} \, ; \, \mathbb{R})$ such that there exist positive constants $C_8 $ and $C_{9,f}$ satisfying with
\begin{equation}
\inf_{|x| \le C_8} |f^\prime (x)|^2 \ge C_{9,f}. 
\end{equation}
Then, we can find the inverse function of $f$ around the origin. 
Define 
$$
Y(t) := \frac{ \ 1 \ }{t} \int_0^t f(X(s)) \, \mbox{d}s \quad (t \in (0,T]). 
$$

\begin{Thm}\label{Thm:smooth density for Y(t)}
Suppose the uniformly elliptic condition \eqref{UE condition} on $A_1$ stated in Lemma \ref{Lem:density for X(t)}. 
Then, for each $0 < t \le T$, the probability law of $Y(t)$ admits a smooth density with respect to the Lebesgue measure on $\mathbb{R}$.  
\end{Thm}
{\it{Proof}}. 
We have only to check the negative-order integrability of the Malliavin covariance $\tilde{V}(t)$ of $Y(t)$, because the Malliavin smoothness of $Y(t)$ can be derived by the one of $X(s)$ for each $0 \le s \le t$. 
See \cite{Nualart}. 

Firstly, we shall compute the Malliavin covariance $\tilde{V}(t) := \langle DY(t), DY(t) \rangle_{\mathbb{H}_0^1}$.  
Since 
\begin{align*}
\frac{ \mbox{d}}{\mbox{d}u} D_u Y(t) 
& = \frac{ \ 1 \ }{t} \, \int_0^t f^\prime (X(s)) \, \frac{ \mbox{d}}{\mbox{d}u} D_u X(s) \, \mbox{d}s \\
& = \frac{ \ 1 \ }{t} \, \int_u^t f^\prime (X(s)) \, Z(s,u) \, A_1 (u, X_u) \, \mbox{d}s, 
\end{align*}
for $0 \le u \le t$, the Malliavin covariance $\tilde{V}(t)$ of $Y(t)$ can be computed as follows: 
\begin{align*}
\tilde{V}(t) 
& = \int_0^t \left( \frac{\mbox{d}}{\mbox{d}u} D_u Y(t) \right)^2 \, \mbox{d}u \\
& = \int_0^t \left\{ \frac{ \ 1 \ }{t} \left( \int_u^t f^\prime (X(s)) \, Z(s,u) \, \mbox{d}s \right) \, A_1 (u, X_u) \right\}^2 \mbox{d}u.  
\end{align*}

Secondly, we shall check the negative-order integrability of $\tilde{V}(t)$. 
In order to do it, it is sufficient to observe $\mathbb{E} \big[ \exp (- \lambda \, \tilde{V}(t) \big] = o(\lambda^{-p})$ as $\lambda \to +\infty$ for any $p > 1$, because 
$$
\mathbb{E} \big[ \tilde{V}(t)^{-p} \big] 
= \frac{1}{\Gamma (p)} \, \int_0^{+\infty} \lambda^{p-1} \, \mathbb{E} \big[ \exp (- \lambda \, \tilde{V}(t) ) \big] \, \mbox{d} \lambda. 
$$
Let $0 < \gamma < 1/3$ be a constant, and $\lambda > 1$ sufficiently large. 
Write $t_\lambda := t - \lambda^{- \gamma}$ and 
$$
\Omega_3 := \left\{ \sup_{t_\lambda \le s \le t} |X(s)| \le C_8 \right\} 
\cap \left\{ \sup_{t_\lambda \le s \le t} \big| Z \big( s, \theta ( t_\lambda, t ) \big) -1 \big|^2 \le \frac{ \ 1 \ }{4} \right\}. 
$$
We shall remark that 
\begin{align*}
\mathbb{P} [ \Omega_3^c ]
& \le \mathbb{P} \left[ \sup_{t_\lambda \le s \le t} |X(s)| > C_8 \right] 
+ \mathbb{P} \left[ \sup_{t_\lambda \le s \le t} \big| Z \big( s, \theta ( t_\lambda, t ) \big) - 1 \big|^2 >\frac{ \ 1 \ }{4} \right] \\
& \le C_8^{-p} \, \mathbb{E} \left[ \sup_{t_\lambda \le s \le t} |X(s)|^p \right] 
+ 4^p \, \mathbb{E} \left[ \sup_{t_\lambda \le s \le t} \big| Z \big( s, \theta ( t_\lambda, t ) \big) - 1 \big|^{2p} \right] \\
& \le C_{1,p, \eta, T} \, C_8^{-p} \, \lambda^{- \gamma p /2} 
+ C_{4, p, T} \, 4^p \, \lambda^{- \gamma p}
\end{align*}
from the Chebyshev inequality. 
Moreover, the mean value theorem tells us to see that 
\begin{align*}
\tilde{V}(t) 
& \ge \frac{ \ 1 \ }{t^2} \int_{t_\lambda}^t \left\{ \left( \int_u^t f^\prime (X(s)) \, Z(s,u) \, \mbox{d}s \right) \, A_1 (u, X_u) \right\}^2 \mbox{d}u \\
& = \frac{(1-\theta)^2}{t^2} \, \lambda^{- 3 \gamma} \, \left\{ f^\prime (X(s)) \, Z(s,u) \, A_1 (u, X_u) \right\}^2  \big|_{u = \theta (t_\lambda, t), s = \delta (u, t)} \\
& \ge \frac{(1-\theta)^2}{t^2} \, \lambda^{- 3 \gamma}  \, \left( \inf_{|x| \le C_{11}} |f^\prime (x)|^2 \right) \\
& \qquad \times \left( \inf_{0 \le u \le T} \inf_{g \in C([-r,0] \, ; \, \mathbb{R})} |A_1 (u, g)|^2 \right) \ Z \Big( \delta \big( \theta ( t_\lambda, t),t \big), \theta ( t_\lambda, t) \Big)^2 \\
& \ge \frac{(1-\theta)^2}{t^2} \, \lambda^{- 3 \gamma}  \, C_{9,f} \, C_2 \, \left( \frac{ \ 1 \ }{2} - \sup_{ \theta ( t_\lambda, t) \le s \le t} \big| Z \big( s, \theta ( t_\lambda, t ) \big) - 1 \big|^2 \right) \\
& \ge \frac{(1-\theta)^2}{4 \, t^2} \, C_{9,f} \, C_2 \, \lambda^{- 3 \gamma} 
\end{align*}
on $\Omega_3$, where $0 < \theta < 1$ and $0 < \delta < 1$ are constants, $\theta(t_\lambda,t) := \theta t + (1-\theta) t_\lambda$, and $\delta (u,t) := \delta t + (1-\delta) u$. 
In the fifth inequality, we have used the assumption on the function $f$, and the uniformly elliptic condition \eqref{UE condition} on $A_1$. 
Write $\hat{\mathbb{E}} \big[ \ \cdot \ \big] = \mathbb{E} \left[ \ \cdot \ \mathbb{I}_{\Omega_3} \right]$. 
Then, we can get 
\begin{align*}
\mathbb{E} \big[ \exp (- \lambda \, \tilde{V}(t)) \big] 
& \le \hat{\mathbb{E}} \big[ \exp (- \lambda \, \tilde{V}(t)) \big] 
+ \mathbb{P} [ \Omega_3^c ] \\
& \le \exp \left[ - \frac{(1-\theta)^2}{4 t^2} \, C_{9,f} \, C_2 \, \lambda^{1 - 3 \gamma} \right] 
+ C_{1, p, \eta, T} \, C_8^{-p} \, \lambda^{- \gamma p /2} 
+ C_{4, p, T} \, 4^p \, \lambda^{- \gamma p} \\
& = o(\lambda^{- C_{10} \, \, p})
\end{align*}
as $\lambda \to +\infty$ for any $p > 1$, which is our desired conclusion. 
\hfill $\square$

\begin{Rem}
{\rm{
We can also obtain Theorem \ref{Thm:smooth density for Y(t)} from the viewpoint of Theorem \ref{Thm:n-state}. 
Let $0 = t_0 < t_1 < \cdots < t_N = t$ such that $\| \Delta \| := \max_{1 \le k \le N} (t_{k} - t_{k-1})$ tends to $0$ as $N \to +\infty$. 
Write 
$$
Y(t) = \frac{ \ 1 \ }{t} \int_0^t f(X(s)) \, \mbox{d}s, \quad 
Y_N(t_1, \dots, t_N) := \frac{ \ 1 \ }{N} \sum_{k=1}^N f(X(t_k)). 
$$
The corresponding Malliavin covariances to $Y(t)$ and $Y_N(t_1, \dots, t_N)$ are 
\begin{align*}
& \tilde{V} (t) = \int_0^T \left\{ \frac{ \ 1 \ }{t} \int_0^t f^\prime (X(s)) \, Z(s,u) \, A_1 (u, X_u) \, \mathbb{I}_{(u \le s)} \, \mbox{d}s \right\}^2 \mbox{d}u, \\
& \tilde{V}_N (t_1, \dots, t_N) = \int_0^T \left\{ \frac{ \ 1 \ }{N} \sum_{k=1}^N f^\prime (X(t_k)) \, Z(t_k, u) \, A_1 (u, X_u) \, \mathbb{I}_{(u \le t_k)} \right\}^2 \mbox{d}u. 
\end{align*}
As pointed out in the proof of Theorem \ref{Thm:smooth density for Y(t)}, the Malliavin smoothness of $Y(t)$ and $Y_N(t_1, \dots, t_N)$ can be checked directly from the one of $X(s)$ for $0 \le s \le t$. 
On the other hand, in order to study the negative-order moment of $\tilde{V}_N (t)$, it is sufficient to give the one ov $\tilde{V}_N (t_1, \dots, t_N)$, because the dominated convergence theorem and the Fatou lemma lead us to see that   
\begin{align*}
\mathbb{E} \big[ \tilde{V}(t)^{-q} \big] 
& = \frac{1}{\Gamma(q)} \, \int_0^{+\infty} \lambda^{q-1} \, \mathbb{E} \big[ \exp (- \lambda \, \tilde{V}(t)) \big] \, \mbox{d} \lambda \\
& = \frac{1}{\Gamma(q)} \, \int_0^{+\infty} \lambda^{q-1} \, \lim_{N \to +\infty} \mathbb{E} \big[ \exp (- \lambda \, \tilde{V}_N(t_1, \dots, t_N)) \big] \, \mbox{d} \lambda \\
& \le \liminf_{N \to +\infty} \frac{1}{\Gamma(q)} \, \int_0^{+\infty} \lambda^{q-1} \, \mathbb{E} \big[ \exp (- \lambda \, \tilde{V}_N(t_1, \dots, t_N)) \big] \, \mbox{d} \lambda \\
& \le \liminf_{N \to +\infty} \, \mathbb{E} \big[ \tilde{V}_N(t_1, \dots, t_N)^{-q} \big].  
\end{align*}

Now, we shall study the estimate of $\mathbb{E} \big[ \tilde{V}_N(t_1, \dots, t_N)^{-q} \big]$. 
Let $q > 1$ be arbitrary and $0 < \sigma < 1/3$ a constant. 
Write $N_{\sigma, t} := \big( N / t \big)^{1/\sigma}$. 
Remark that 
\begin{align*}
\mathbb{E} \big[ \tilde{V}_N(t_1, \dots, t_N)^{-q} \big] 
& = \frac{1}{\Gamma (q)} \int_0^{+\infty} \lambda^{q-1} \, \mathbb{E} \big[ \exp (- \lambda \tilde{V}_N(t_1, \dots, t_N)) \big] \, \mbox{d} \lambda \\
& \le \frac{1}{\Gamma (q)} \int_0^1 \lambda^{q-1} \, \mbox{d} \lambda 
+ \frac{1}{\Gamma (q)} \int_1^{N_{\sigma,t}} \lambda^{q-1} \, \mathbb{E} \big[ \exp ( - \lambda \tilde{V}_N(t_1, \dots, t_N)) \big] \, \mbox{d} \lambda \\
& \qquad + \frac{1}{\Gamma (q)} \int_{N_{\sigma,t} }^{+\infty} \lambda^{q-1} \, \mathbb{E} \big[ \exp ( - \lambda \tilde{V}_N(t_1, \dots, t_N)) \big] \, \mbox{d} \lambda \\
& = : I_{1,N} + I_{2,N} + I_{3,N}. 
\end{align*}
First of all, the estimate of $I_{1,N}$ is trivial: 
$$
I_{1,N} 
\equiv \frac{1}{\Gamma (q)} \int_0^1 \lambda^{q-1} \, \mbox{d} \lambda 
= \frac{ \ 1 \ }{\Gamma (q + 1)} 
< +\infty. 
$$
Secondly, we shall consider the estimate of $I_{2,N}$. 
Since $Z(t,u) = 0$ for $t < u$, the mean value theorem implies that 
\begin{align*}
\tilde{V}_N(t_1, \dots, t_N) 
& = \int_0^{t} \left\{ \frac{ \ 1 \ }{N} \sum_{k=1}^N f^\prime (X(t_k)) \, Z(t_k, u) \, A_1 (u, X_u) \, \mathbb{I}_{(u \le t_k)} \right\}^2 \mbox{d}u \\
& \ge \int_{t_{N-1}}^{t} \left\{ \frac{ \ 1 \ }{N} \sum_{k=1}^N f^\prime (X(t_k)) \, Z(t_k, u) \, A_1 (u, X_u) \, \mathbb{I}_{(u \le t_k)} \right\}^2 \mbox{d}u \\
& = \int_{t_{N-1}}^{t} \left\{ \frac{ \ 1 \ }{N} f^\prime (X(t)) \, Z(t, u) \, A_1 (u, X_u) \right\}^2 \mbox{d}u \\
& = \left\{ \frac{ \ 1 \ }{N} f^\prime \big( X(t) \big) \, Z \big( t, \theta_2 (t_{N-1}, t) \big) \, A_1 \big( \theta_2 (t_{N-1}, t), X_{\theta_2 (t_{N-1}, t)} \big) \right\}^2 \, N_{\sigma,t}^{-\sigma}, 
\end{align*}
where $0 < \theta_2 < 1$ is a constant, and $\theta_2 (s,t) := s + \theta_2 \, (t-s)$. 
Denote by 
$$
\Omega_4 := \left\{ \sup_{t_{N-1} \le s \le t} \big| X(s) \big| \le C_8 \right\} \cap \left\{ \sup_{\theta_2 (t_{N-1}, t) \le s \le t} \big| Z \big( s, \theta_2 (t_{N-1}, t) \big) - 1 \big|^2 \le \frac{ \ 1 \ }{4} \right\}. 
$$
We shall remark that 
\begin{align*}
\mathbb{P} [ \Omega_4^c ]
& \le \mathbb{P} \left[ \sup_{t_{N-1} \le s \le t} |X(s)| > C_8 \right] 
+ \mathbb{P} \left[ \sup_{\theta_2 (t_{N-1}, t) \le s \le t} \big| Z \big( s, \theta_2 (t_{N-1}, t) \big) - 1 \big|^2 >\frac{ \ 1 \ }{4} \right] \\
& \le C_8^{-p} \, \mathbb{E} \left[ \sup_{t_{N-1} \le s \le t} |X(s)|^p \right] 
+ 4^p \, \mathbb{E} \left[ \sup_{\theta_2 (t_{N-1}, t) \le s \le t} \big| Z \big( s, \theta_2 ( t_{N-1}, t ) \big) - 1 \big|^{2p} \right] \\
& \le C_{1,p, \eta, T} \, C_8^{-p} \, N_{\sigma,t}^{- \sigma p /2}
+ C_{4, p, T} \, 4^p \, N_{\sigma,t}^{-\sigma p}
\end{align*}
from the Chebyshev inequality. 
Moreover, we see that, on $\Omega_4$, 
\begin{align*}
\tilde{V}_N (t_1, \dots, t_N) 
& \ge \left\{ \frac{ \ 1 \ }{N} f^\prime \big( X(t) \big) \, Z \big( t, \theta_2 (t_{N-1}, t) \big) \, A_1 \big( \theta_2 (t_{N-1}, t), X_{\theta_2 (t_{N-1}, t)} \big) \right\}^2 \, N_{\sigma,t}^{-\sigma} \\
& \ge \frac{C_{9, f} \, C_2 \ N_{\sigma,t}^{-\sigma}}{N^2} \, Z \big( s, \theta_2 (t_{N-1}, t) \big)^2 \\
& \ge \frac{C_{9, f} \, C_2 \ N_{\sigma,t}^{-\sigma}}{4 N^2}
\end{align*}
from the assumption on the function $f$ and the condition \eqref{UE condition}. 
Write $\hat{E} \big[ \ \cdot \ \big] := \mathbb{E} \big[ \ \cdot \ \mathbb{I}_{\Omega_4} \big]$. 
Hence, it holds that 
\begin{align*}
\mathbb{E} \big[ \exp (- \lambda \, \tilde{V}_N(t_1, \dots, t_N)) \big] 
& \le \hat{\mathbb{E}} \big[ \exp (- \lambda \, \tilde{V}_N(t_1, \dots, t_N)) \big] 
+ \mathbb{P} \big[ \Omega_4^c \big] \\
& \le \exp \left( - \frac{C_{9, f} \, C_2 \ N_{\sigma,t}^{-\sigma}}{4 N^2} \, \lambda \right) 
+ C_{1, p, \eta, T} \, C_8^{-p} \, N_{\sigma,t}^{- \sigma p /2} 
+ C_{4, p, T} \, 4^p \, N_{\sigma,t}^{- \sigma p} \\
& \le \exp \left( - \frac{C_{9, f} \, C_2 \ N_{\sigma,t}^{-\sigma}}{4 N^2} \, \lambda \right) 
+ \left( C_{1, p, \eta, T} \, C_8^{-p} + C_{4, p, T} \, 4^p \right) N_{\sigma,t}^{- \sigma p /2}. 
\end{align*}
Let $p > \big( 2q / \sigma \big) \vee 1$ be arbitrary. 
Thus, the mean value theorem enables us to get 
\begin{align*}
I_{2,N} 
& \equiv \frac{1}{\Gamma (q)} \int_1^{N_{\sigma,t} } \lambda^{q-1} \, \mathbb{E} \big[ \exp ( - \lambda \tilde{V}_N(t_1, \dots, t_N)) \big] \, \mbox{d} \lambda \\
& \le \frac{1}{\Gamma (q)} \int_1^{N_{\sigma,t} } \lambda^{q-1} \, \left\{ \exp \left( - \frac{C_{9, f} \, C_2 \ N_{\sigma,t}^{-\sigma}}{4 N^2} \, \lambda \right) 
+ \left( C_{1, p, \eta, T} \, C_8^{-p} + C_{4, p, T} \, 4^p \right) N_{\sigma,t}^{- \sigma p /2} \right\} \mbox{d} \lambda \\
& \le \frac{\lambda^{q-1}}{\Gamma (q)} \, \exp \left( - \frac{C_{9, f} \, C_2 \ N_{\sigma,t}^{-\sigma}}{4 N^2} \, \lambda \right) \bigg|_{ \lambda = (1-\delta_2) + \delta_2 \, N_{\sigma,t}  } \, \big( N_{\sigma,t}  -1 \big) \\
& \qquad + \frac{1}{\Gamma (q)} \int_1^{N_{\sigma,t} } \left( C_{1, p, \eta, T} \, C_8^{-p} + C_{4, p, T} \, 4^p \right) \lambda^{q - \sigma p/2 -1} \, \mbox{d} \lambda \\
& \le \frac{C_{11,q}}{\Gamma (q)} \, \big( 1 + N_{\sigma,t} ^q \big) \, \exp \left[ - \frac{C_{9, f} \, C_2 \, \delta_2 \, N_{\sigma,t}^{1-\sigma}}{4 N^2} \right] \\
& \qquad + \frac{1}{\Gamma (q)} \int_1^{+\infty} \left( C_{1, p, \eta, T} \, C_8^{-p} 
+ C_{4, p, T} \, 4^p \right) \lambda^{q - \sigma p/2 -1} \, \mbox{d} \lambda \\
& = \frac{C_{11, q}}{\Gamma (q)} \, \big( 1 + N_{\sigma,t} ^q \big) \, \exp \left[ - \frac{C_{9, f} \, C_2 \, \delta_2 \, N_{\sigma,t}^{1-\sigma}}{4 N^2} \right] 
+ \frac{C_{1, p, \eta, T} \, C_9^{-p} + C_{4, p, T} \, 4^p}{( \sigma p / 2 - q) \, \Gamma (q)} \\
& \to \frac{C_{1, p, \eta, T} \, C_8^{-p} + C_{4, p, T} \, 4^p}{( \sigma p / 2 - q) \, \Gamma (q)} \quad (N \to +\infty), 
\end{align*}
because of $0 < \sigma < 1/3$, where $0 < \delta_2 < 1$ is a constant. 

Thirdly, we shall consider the estimate of $I_{3,N}$. 
Write $t_\lambda := t - \lambda^{- \sigma}$, and let $\lambda \ge N_{\sigma, t}$ be sufficiently large. 
The mean value theorem implies that 
\begin{align*}
\tilde{V}_N(t_1, \dots, t_N) 
& = \int_0^t \left\{ \frac{ \ 1 \ }{N} \sum_{k=1}^N f^\prime (X(t_k)) \, Z(t_k, u) \, A_1 (u, X_u) \, \mathbb{I}_{(u \le t_k)} \right\}^2 \mbox{d}u \\
& \ge \int_{t_\lambda}^t \left\{ \frac{ \ 1 \ }{N} \sum_{k=1}^N f^\prime (X(t_k)) \, Z(t_k, u) \, A_1 (u, X_u) \, \mathbb{I}_{(u \le t_k)} \right\}^2 \mbox{d}u \\
& = \int_{t_\lambda}^t \left\{ \frac{ \ 1 \ }{N} f^\prime (X(t)) \, Z(t, u) \, A_1 (u, X_u) \right\}^2 \mbox{d}u \\
& = \left\{ \frac{ \ 1 \ }{N} f^\prime \big( X(t) \big) \, Z \big( t, \theta_3(t_\lambda, t) \big) \, A_1 \big( \theta_3 (t_\lambda, t), X_{\theta_3 (t_\lambda, t)} \big) \right\}^2 \lambda^{-\sigma}, 
\end{align*}
where $0 < \theta_3 < 1$ is a constant, and $\theta_3 (s,t) := s + \theta_3 \, (t-s)$. 
Denote by 
\begin{align*}
\Omega_{5,\lambda} &:= \left\{ \sup_{t_\lambda \le s \le t} \big| X(s) \big| \le C_8 \right\} \cap \left\{ \sup_{\theta_3 (t_\lambda, t) \le s \le t} \big| Z \big( s, \theta_3 (t_\lambda, t) \big) - 1 \big|^2 \le \frac{ \ 1 \ }{4} \right\}. 
\end{align*}
Then, we can get 
\begin{align*}
\tilde{V}_N (t_1, \dots, t_N) 
& \ge \left\{ \frac{ \ 1 \ }{N} f^\prime \big( X(t) \big) \, Z \big( t, \theta_3 (t_\lambda, t) \big) \, A_1 \big( \theta_3 (t_\lambda, t), X_{\theta_3 (t_\lambda, t)} \big) \right\}^2 \lambda^{-\sigma} \\
& \ge \frac{C_{9, f} \, C_2}{N^2} \, Z \big( s, \theta_3 (t_\lambda, t) \big)^2 \, \lambda^{-\sigma} \\
& \ge \frac{C_{9, f} \, C_2}{4 t^2} \, \lambda^{-3\sigma}
\end{align*}
on $\Omega_{5, \lambda}$, from the assumption on the function $f$ and the uniformly elliptic condition \eqref{UE condition}, because of $t/N = N_{\sigma, t}^{-\sigma} \ge \lambda^{-\sigma}$. 
Define $\hat{\mathbb{E}}_\lambda \big[ \ \cdot \ \big] := \mathbb{E} \big[ \ \cdot \ \mathbb{I}_{\Omega_{5,\lambda}} \big]$. 
Then, we have 
\begin{align*}
\mathbb{E} \big[ \exp ( - \lambda \tilde{V}_N(t_1, \dots, t_N)) \big] 
& \le \hat{\mathbb{E}}_\lambda \big[ \exp ( - \lambda \tilde{V}_N(t_1, \dots, t_N)) \big] 
+ \mathbb{P} \big[ \Omega_{5,\lambda}^c \big] \\
& \le \exp \left( - \frac{C_{9, f} \, C_2}{4 t^2} \, \lambda^{1-3\sigma} \right) 
+ C_{1, p,\eta,T} \, C_8^{-p} \, \lambda^{- \sigma p /2} 
+ C_{4, p,T} \, 4^p \, \lambda^{-\sigma p} \\
& \le  \exp \left( - \frac{C_{9, f} \, C_2}{4 t^2} \, \lambda^{1-3\sigma} \right) 
+ \left( C_{1, p,\eta,T} \, C_8^{-p} + C_{4, p,T} \, 4^p \right) \lambda^{-\sigma p/2}
\end{align*}
from the Chebyshev inequality. 
Hence, it holds that 
\begin{align*}
I_{3,N} 
& \equiv \frac{1}{\Gamma (q)} \int_{t_{\sigma, N}}^{+\infty} \lambda^{q-1} \, \mathbb{E} \big[ \exp ( - \lambda \tilde{V}_N(t_1, \dots, t_N)) \big] \, \mbox{d} \lambda \\
& \le \frac{1}{\Gamma (q)} \int_{t_{\sigma, N}}^{+\infty} \lambda^{q-1} \, \left\{ \exp \left( - \frac{C_{9, f} \, C_2}{4 t^2} \, \lambda^{1-3\sigma} \right) 
+ \left( C_{1, p,\eta,T} \, C_8^{-p} + C_{4, p,T} \, 4^p \right) \lambda^{-\sigma p/2} \right\} \, \mbox{d} \lambda \\
& \le \frac{1}{\Gamma (q)} \int_1^{+\infty} \lambda^{q-1} \, \left\{ \exp \left( - \frac{C_{9, f} \, C_2}{4 t^2} \, \lambda^{1-3\sigma} \right) 
+ \left( C_{1, p,\eta,T} \, C_8^{-p} + C_{4, p,T} \, 4^p \right) \lambda^{-\sigma p/2} \right\} \, \mbox{d} \lambda \\
& \le \frac{( 4t^2 )^{q / (1-3\sigma)} \, \Gamma (q / (1-3\sigma))}{(1 - 3 \sigma) \, (C_{9, f} \, C_2)^{q / (1-3\sigma)} \Gamma (q)}
+ \frac{C_{1, p,\eta,T} \, C_8^{-p} + C_{4, p,T} \, 4^p}{(q - \sigma p /2) \, \Gamma (q)}. 
\end{align*}

Therefore, we see that 
\begin{align*}
\mathbb{E} \big[ \tilde{V}(t)^{-q} \big] 
& \le \liminf_{N \to +\infty} \big( I_{1,N} + I_{2,N} + I_{3,N} \big) \\
& \le \frac{ \ 1 \ }{\Gamma (q+1)} 
+ \frac{C_{1, p, \eta, T} \, C_8^{-p} + C_{4, p, T} \, 4^p}{( \sigma p / 2 - q) \, \Gamma (q)} \\
& \quad + \left\{ \frac{( 4t^2 )^{q / (1-3\sigma)} \, \Gamma (q / (1-3\sigma))}{(1 - 3 \sigma) \, (C_{9, f} \, C_2)^{q / (1-3\sigma)} \Gamma (q)}
+ \frac{C_{1, p,\eta,T} \, C_8^{-p} + C_{4, p,T} \, 4^p}{(q - \alpha p /2) \, \Gamma (q)} \right\} 
< +\infty, 
\end{align*}
which is our desired one. 
\hfill $\square$
}}
\end{Rem}

\section{Delayed Black-Scholes models and sensitivity analysis}\label{Section 5}

In this section, we shall apply our studies to the option pricing of the asset price dynamics model with delayed effects. 
See \cite{Arriojas-Hu-Mohammed-Pap, Chang-Pang-Pemy, Chang-Youree-1, Chang-Youree-2, MS} on details. 
Let $A_i \in C_{0+,b}^1 ( \mathbb{R} \, ; \, \mathbb{R} ) \ (i = 0, \, 1)$ with the uniformly elliptic condition of the form: there exists a positive constant $C_{12}$ with 
\begin{equation}\label{condition:UE}
\inf_{y \in \mathbb{R}} A_1 ( y )^2 \ge C_{12}. 
\end{equation}
Let $x > 0$ be a constant. 
Consider the $\mathbb{R}$-valued process $X = \big\{ X(t) \, ; \, -r \le t \le T \big\}$ determined by the stochastic delay differential equation of the form: 
\begin{equation}\label{eq:SDDE X}
X(t) = 
\begin{cases}
x & (-r \le t \le 0), \\
\displaystyle x + \int_0^t A_0 \big( X(s-r) \big) \, X(s) \, \mbox{d}s 
+ \int_0^t A_1 \big( X(s-r) \big) \, X(s) \, \mbox{d}W(s) & (0 < t \le T). 
\end{cases}
\end{equation}
The existence of the unique solution to \eqref{eq:SDDE X} can be checked easily. 
In fact, it is trivial on the interval $[-r,0]$. 
Since 
$$
X(t) = x + \int_0^t A_0 (x ) \, X(s) \, \mbox{d}s 
+ \int_0^t A_1 ( x ) \, X(s) \, \mbox{d}W(s)
$$
for $0 < t \le r$, we can derive our conclusion from the Lipschitz condition and the linear growth condition on the coefficients. 
Iterating such argument enables us to get our desired assertion on each intervals $[kr,(k+1)r] \ (k \in \mathbb{N} \cup \{ 0 \})$. 
Moreover, the equation \eqref{eq:SDDE X} can be solved as 
\begin{equation}
X(t) = 
\begin{cases}
x \qquad \qquad \qquad \qquad \qquad \qquad \qquad \qquad \qquad \qquad \qquad \qquad \qquad \qquad (-r \le t \le 0), \\
\displaystyle x \, \exp \! \Bigg[ \! \int_0^t \!\! A_1 \big( X(s-r) \big) \mbox{d}W(s) 
\! + \! \int_0^t \! \Bigg\{ \! A_0 \big( X(s-r) \big) - \frac{A_1 \big( X(s-r) \big)^2}{2} \! \Bigg\} \mbox{d}s \Bigg] \ (0 \le t \le T), 
\end{cases}
\end{equation}
which implies that $X(t) > 0$ a.s. 
The process $X = \big\{ X(t) \, ; \, -r \le t \le T \big\}$ is called {\it{the delayed Black-Scholes model}} in mathematical finance (cf. \cite{Arriojas-Hu-Mohammed-Pap, Chang-Pang-Pemy, Chang-Youree-1, Chang-Youree-2, MS}). 
\begin{Rem}\label{Rem:density of SDDE X(t)}
{\rm{
Although the coefficients of the equation \eqref{eq:SDDE X} do not satisfy the uniformly elliptic condition \eqref{UE condition} as stated in Lemma \ref{Lem:density for X(t)}, the probability law of $X(t)$ admits a smooth density with respect to the Lebesgue measure on $\mathbb{R}$. 
Denote by $\tilde{x} = \log x$ and $\tilde{A}_i (\tilde{y}) = A_i (e^{\tilde{y}}) \ (i = 0, \, 1)$. 
Then, the process $\tilde{X} = \big\{ \tilde{X}(t) := \log X(t) \, ; \, -r \le t \le T \big\}$ satisfies 
\begin{equation}\label{eq:SDDE of tilde{X}}
\tilde{X}(t) = 
\begin{cases}
\tilde{x} & (-r \le t \le 0), \\
\displaystyle \tilde{x} + \int_0^t \tilde{A}_1 (\tilde{X}(s-r)) \, \mbox{d}W(s) 
+ \int_0^t \left\{\tilde{A}_0 (\tilde{X}(s-r)) - \frac{A_1 (\tilde{X}(s-r))}{2} \right\} \, \mbox{d}s & (0 \le t \le T). 
\end{cases}
\end{equation}
Since the coefficient $\tilde{A}_1$ satisfies 
$$
\inf_{\tilde{y} \in \mathbb{R}} \tilde{A}_1 (\tilde{y})^2 
= \inf_{\tilde{y} \in \mathbb{R}} A_1 (e^{\tilde{y}})^2 
\ge C_{12}
$$
from \eqref{condition:UE}, we can conclude from Lemma \ref{Lem:density for X(t)} that the probability law of the $\mathbb{R}$-valued random variable $\tilde{X}(t)$ admits a smooth density $\tilde{p}_t (\tilde{y})$ with respect to the Lebesgue measure on $\mathbb{R}$. 
Hence, the density of the probability law of $X(t) = \exp (\tilde{X}(t))$ is 
$$
p_t (y) = \frac{\tilde{p}_t (\log y)}{y} \,  \quad (y > 0), 
$$
which is smooth in $y > 0$. 
\hfill $\square$
}}
\end{Rem}

Let $R > 0$ be a constant, which denotes the rate of return of a riskless asset. 
Denote by $B = \big\{ B(t) \, ; \, -r \le t \le T \big\}$ the riskless asset price process, which is given by 
$$
B(t) = \mathbb{I}_{[-r,0]} (t) + e^{Rt} \, \mathbb{I}_{(0,T]}(t). 
$$
Write $\bar{X} = \big\{ \bar{X}(t):= X(t) / B(t) \, ; \, -r \le t \le T \big\}$, which is called {\it{the discounted stock price process}}. 
Then, the It\^o formula leads us to see that, for $0 < t \le T$, 
\begin{align*}
\mbox{d} \bar{X}(t) 
& = \left\{ A_0 \big( X(t-r) \big) - R \right\} \, \bar{X}(t) \, \mbox{d}t 
+ A_1 \big( X(t-r) \big)  \, \bar{X}(t) \, \mbox{d}W(t) \\
& = A_1 \big( X(t-r) \big)  \, \bar{X}(t) \, \big\{ \mbox{d}W(t) - \Sigma (t) \, \mbox{d}t \big\}, 
\end{align*}
where $\Sigma(t) = - \big\{ A_0 \big( X(t-r) \big)  - R \big\} / A_1 \big( X(t-r) \big) $. 
Define the process $M = \big\{ M(t) \, ; \, 0 \le t \le T \big\}$ by 
$$
M(t) := \exp \left[ \int_0^t \Sigma(s) \, \mbox{d} W(s) - \frac{ \ 1 \ }{2} \int_0^t \Sigma(s)^2 \, \mbox{d}s \right]. 
$$
Then, the process $M$ is a square-integrable $(\mathcal{F}_t)$-martingale, because 
\begin{align*}
\mathbb{E} \big[ M(t)^2 \big] 
& \le \mathbb{E} \left[ \exp \left\{ 8 \int_0^t \Sigma(s)^2 \, \mbox{d}s \right\} \right]^{1/2} 
< +\infty 
\end{align*}
from the boundedness of $A_0$ and the uniformly elliptic condition \eqref{condition:UE} on $A_1$. 
In particular, we have $\mathbb{E} [ M(T)] = 1$. 
Then, the measure $\mbox{d} \mathbb{Q} := M(T) \, \mbox{d} \mathbb{P}$ is also the probability one on the measurable space $(\Omega, \mathcal{F})$, and the Girsanov theorem tells us to see that the process $\tilde{W} = \big\{ \tilde{W}(t) := W(t) - \int_0^t \Sigma (s) \, \mbox{d}s \, ; \, 0 \le t \le T \big\}$ is also a Brownian motion starting from the origin under the measure $\mbox{d} \mathbb{Q}$.  
Let $Z$ be a $\mathcal{F}_T$-measurable, non-negative and integrable random variable, which is called a contingent claim on the process $X = \big\{ X(t) \, ; \, -r \le t \le T \big\}$. 


\begin{Prop}\label{Prop:completeness of market}
The market $(X,B) = \big\{ (X(t), B(t)) \, ; \, 0 \le t \le T \big\}$ is complete. 
\end{Prop}
{\it{Proof}}. 
Remark that $\mathcal{F}_t^{\tilde{W}} = \mathcal{F}_t$. 
Since the process $L = \big\{ L(t)  := \mathbb{E}_{\mathbb{Q}} \big[ e^{-RT} \, Z | \mathcal{F}_t \big] \, ; \, 0 \le t \le T \big\}$  is an $(\mathcal{F}_t)$-martingale under $\mbox{d} \mathbb{Q}$, we can find an $(\mathcal{F}_t)$-predictable, square-integrable process $\varphi = \big\{ \varphi (t) \, ; \, 0 \le t \le T \big\}$ such that 
$$
L(t) = L(0) + \int_0^t \varphi (s) \, \mbox{d} \tilde{W}(s)
$$
by the martingale representation theorem (cf. \cite{Ikeda-Watanabe}-Theorem II-6.6, p.80). 
Set 
\begin{gather*}
\pi_X(t) := \frac{ \varphi (t)}{A_1 \big( X(t-r) \big)  \, \tilde{X}(t)}, \ \ 
\pi_B (t) := L(t) - \pi_X(t) \, \tilde{X}(t), \ \ 
V(t) := \pi_B(t) \, B(t) + \pi_X(t) \, X(t). 
\end{gather*}
Then, since $V(t) = B(t) \, L(t)$, the It\^o formula implies that 
\begin{align*}
\mbox{d}V(t) 
& = B(t)\, \mbox{d}L(t) + L(t) \, \mbox{d} B(t) \\
& = B(t) \, \varphi (t) \, \mbox{d} \tilde{W}(t) + \big\{\pi_B (t) + \pi_X(t) \, \tilde{X}(t) \big\} \mbox{d}B(t) \\
& = \pi_X (t) \left\{ A_1 \big( X(t-r) \big) \, X(t) \, \big( \mbox{d}W(t) - \Sigma(t) \, \mbox{d}t \big) + R \, X(t) \, \mbox{d}t \right\} + \pi_B(t) \, \mbox{d}B(t) \\
& = \pi_X(t) \, \mbox{d} X(t) + \pi_B(t) \, \mbox{d}B(t), 
\end{align*}
which means that $(\pi_X, \pi_B) = \big\{ (\pi_X (t), \pi_B (t)) \, ; \, 0 \le t \le T \big\}$ is a self-financing strategy. 
Moreover, since 
$$
V(T) = e^{RT} \, L(T) = \mathbb{E}_{\mathbb{Q}} \big[ Z | \mathcal{F}_T \big] = Z, 
$$
the contingent claim $Z$ is attainable. 
Hence, the market $(X,B) = \big\{ (X(t), B(t)) \, ; \, 0 \le t \le T \big\}$ is complete, which is our goal in the proposition. 
\hfill $\square$

\begin{Prop}
For $0 \le t \le T$, it holds that 
\begin{equation}
V(t) = e^{- R(T-t)} \, \mathbb{E}_{\mathbb{Q}} \big[ Z | \mathcal{F}_t \big]. 
\end{equation}
\end{Prop}
{\it{Proof}}. 
As seen in the proof of Proposition \ref{Prop:completeness of market}, it holds that 
$$
\mbox{d} \left( \frac{V(t)}{B(t)} \right) 
= \mbox{d} L(t) 
= \varphi (t) \, \mbox{d} \tilde{W}(t). 
$$
Hence, the discounted process $V / B = \big\{ V(t) / B(t) \, ; \, 0\le t \le T \big\}$ is an $(\mathcal{F}_t)$-martingale under $\mbox{d} \mathbb{Q}$. 
Then, the fair price $V(t)$ of the claim $Z$ is given by 
$$
\frac{V(t)}{B(t)} 
= \mathbb{E}_{\mathbb{Q}} \left[ \frac{V(T)}{B(T)} \bigg| \mathcal{F}_t \right] 
= \mathbb{E}_{\mathbb{Q}} \big[ e^{-RT} \, Z | \mathcal{F}_t \big], 
$$
because of $V(T) = e^{RT} \, L(T) = \mathbb{E}_{\mathbb{Q}} \big[ Z | \mathcal{F}_T \big] = Z$. 
\hfill $\square$

\vspace{5pt}
From now on, we shall consider the case $A_0 \equiv 0$ only, because the Girsanov transform enables us to discuss the general case, if we want. 
For the sake of simplicity of notations, we shall note $\mathbb{P}$ and $W$ instead of $\mathbb{Q}$ and $\tilde{W}$. 
Since 
$$
X(t) = 
\begin{cases}
x & (-r \le t \le 0), \\
\displaystyle x + \int_0^t A_1 \big( X(s-r) \big)  \, X(s) \, \mbox{d}W(s) & (t \in [0,T]), 
\end{cases}
$$
we have 
$$
\partial_x X(t) = 
\begin{cases}
1 & (t \in [-r,0]), \\
\displaystyle 1 + \int_0^t A_1^\prime \big( X(s-r) \big) \, \partial_x X(s-r)  \, X(s) \, \mbox{d}W(s) 
+ \int_0^t A_1 \big( X(s-r) \big) \, \partial_x X(s) \, \mbox{d}W(s) & (0 \le t \le T). 
\end{cases}
$$
Let $U = \big\{ U(t) \, ; \, -r \le t \le T \big\}$ and $\hat{U} = \big\{ \hat{U}(t) \, ; \, -r \le t \le T \big\}$ be $\mathbb{R}$-valued processes determined by the equations: 
\begin{align*}
& U(t) = \hat{U}(t) = 1 \qquad \qquad \qquad \qquad \qquad \qquad \qquad \qquad \qquad \qquad \qquad (-r \le t \le 0), \\
& U(t) = 1 + \int_0^t A_1 \big( X(s-r) \big) \, U(s) \, \mbox{d}W(s) \qquad \qquad \qquad \qquad \qquad \qquad \quad (0 \le t \le T), \\
& \hat{U}(t) = 1 - \int_0^t \hat{U}(s) \, A_1 \big( X(s-r) \big) \, \mbox{d}W(s) 
+ \int_0^t \hat{U}(s) \, A_1 \big( X(s-r) \big)^2 \, \mbox{d}s \qquad (0 \le t \le T). 
\end{align*}
Then, it is clear that $U(t) \, \hat{U}(t) = \hat{U}(t) \, U(t) = 1$, via the It\^o formula, and that 
\begin{align*}
\partial_x X(t) 
& = U(t) \Bigg\{ 
1 + \int_0^t \hat{U}(s) \, A_1^\prime \big( X(s-r) \big) \, \partial_x X(s-r) \, X(s) \, \mbox{d}W(s) \\
& \qquad \qquad - \int_0^t \hat{U}(s) \, A_1 \big( X(s-r) \big) \, A_1^\prime \big( X(s-r) \big) \, \partial_x X(s-r) \, \mbox{d}s \Bigg\} \\
& = : U(t) \, \Lambda(t). 
\end{align*}
Moreover, for $0 \vee (t-r) \le u \le t$, we see that $Z(t,u) = U(t) \, \hat{U}(u)$ is invertible, because 
\begin{align*}
Z(t,u) 
& = 1 + \int_u^t A_1^\prime \big( X(s-r) \big) \, Z (s-r ,u) \, X(s) \, \mbox{d}W(s) 
+ \int_u^t A_1 \big( X(s-r) \big)\, Z(s,u) \, \mbox{d}W(s) \\
& = 1 + \int_u^t A_1 \big( X(s-r) \big) \, Z(s,u) \, \mbox{d}W(s) 
\end{align*}
for $0 \vee (t-r) \le u \le t$, 
$$
U(t) \, \hat{U}(u) 
= 1 + \int_u^t A_1 \big( X(s-r) \big) \, U(s) \, \hat{U}(u) \, \mbox{d}W(s),
$$
and the uniqueness of the solutions. 

Now, we shall state the result on the Greeks computation on the European option with respect to the initial point $x$. 
Denote by $C_{LG} (\mathbb{R} \, ; \, \mathbb{R})$ the set of continuous functions with the linear growth order, and define 
$$
\mathfrak{F} (\mathbb{R} \, ; \, \mathbb{R}) 
= \left\{ \sum_{k=1}^n \alpha_k \, f_k \, \mathbb{I}_{K_k} \, ; \, n \in \mathbb{N}, \,  \alpha_k \in \mathbb{R}, \ f_k \in C_{LG} (\mathbb{R} \, ; \, \mathbb{R}), \ K_k \subset \mathbb{R} \, \mbox{: interval} \right\}. 
$$

\begin{Thm}[Greeks computation of the European option]\label{Thm:European Delta}
For $\Phi \in \mathfrak{F} (\mathbb{R}, \, ; \, \mathbb{R})$, it holds that 
\begin{equation}\label{eq:European Delta}
\partial_x \mathbb{E} \big[ \Phi (X(t)) \big] 
= \mathbb{E} \big[ \Phi (X(t)) \, \Gamma_{E}(t) \big], 
\end{equation}
where 
\begin{equation}
\Gamma_E (t) = \frac{ \ 1 \ }{t \wedge r} \, \delta \left( \frac{U(\cdot) \, \Lambda (t)}{A_1 \big( X (\cdot -r) \big) \, X( \cdot )} \, \mathbb{I}_{[0 \vee (t-r),t]} ( \cdot ) \right), 
\end{equation}
and $\delta ( \cdot )$ is the Skorokhod integral operator. 
\end{Thm}
{\it{Proof}}. 
Consider the case of $\Phi \in C_b^1 (\mathbb{R} \, ; \, \mathbb{R})$. 
Choose $0 \vee (t - r) \le u \le t$. 
Since 
\begin{align*}
\frac{ \mbox{d}}{\mbox{d}u} D_u \big( \Phi(X(t)) \big) 
& = \Phi^\prime (X(t)) \, \frac{ \mbox{d}}{\mbox{d}u} D_u X(t) \\
& = \Phi^\prime (X(t)) \,Z(t,u) \, A_1 \big( X(u-r) \big) \, X(u) \\
& = \Phi^\prime (X(t)) \, U(t) \, \hat{U}(u) \, A_1 \big( X(u-r) \big) \, X(u) \\
& = \partial_x \big( \Phi (X(t)) \big) \, \frac{ \hat{U}(u) \, A_1 \big( X(u-r) \big) \, X(u)}{\Lambda (t)}
\end{align*}
from the chain rule on the operator $D$, the integration by parts formula leads us to get 
\begin{align*}
\partial_x \Big( \mathbb{E} \big[ \Phi (X(t)) \big] \Big) 
& = \mathbb{E} \big[ \partial_x \big( \Phi (X(t)) \big) \big] \\
& = \mathbb{E} \left[ \frac{ \ 1 \ }{t \wedge r} \int_{0 \vee (t-r)}^t \frac{ \mbox{d}}{\mbox{d}u} D_u \big( \Phi (X(t)) \big) \, \frac{U(u) \, \Lambda (t)}{A_1 \big( X(u-r) \big) \, X(u)} \, \mbox{d}u \right] \\
& = \mathbb{E} \left[ \Phi (X(t)) \, \frac{ \ 1 \ }{t \wedge r} \, \delta \left( \frac{U( \cdot ) \, \Lambda (t)}{A_1 \big( X(\cdot -r) \big) \, X( \cdot )} \, \mathbb{I}_{[0 \vee (t-r),t]} ( \cdot ) \right) \right] \\
& = \mathbb{E} \big[ \Phi (X(t)) \, \Gamma_E (t) \big]. 
\end{align*}

In order to extend the class of payoff functions, we have to find a sequence $\big\{ \Phi_n \, ; \, n \in \mathbb{N} \big\}$ such that  
\begin{align*}
& \sup_{x \in K} \left| \mathbb{E} \big[ \Phi_n (X(t)) \big] - \mathbb{E} \big[ \Phi (X(t)) \big] \right| \to 0, \\ 
& \sup_{x \in K} \left| \partial_x \mathbb{E} \big[ \Phi_n (X(t)) \big] - \mathbb{E} \big[ \Phi (X(t)) \, \Gamma_E (t) \big] \right| \to 0
\end{align*}
as $n \to +\infty$, where $K$ is a compact subset in $[0,+\infty)$. 
Hence, it is sufficient to justify 
\begin{equation}\label{L^2-estimate}
\mathbb{E} \big[ | \Phi_n (X(t)) - \Phi (X(t)) |^2 \big] \to 0
\end{equation}
as $n \to +\infty$. 
See \cite{Kawai-Takeuchi, Takeuchi} on details. 

As for $\Phi \in C_b (\mathbb{R} \, ; \, \mathbb{R})$, it is easy to find the sequence $\big\{ \Phi_n \, ; \, n \in \mathbb{N} \big\}$ in $C_b^1 (\mathbb{R} \, ; \, \mathbb{R})$ satisfying with \eqref{L^2-estimate}. 
When $\Phi$ is the indicator function, we can approximate $\Phi$ by a sequence $\big\{ \Phi_n \, ; \, n \in \mathbb{N} \big\}$ in $C_b ([0,+\infty) \, ; \, \mathbb{R})$. 
In order to give the convergence \eqref{L^2-estimate}, we need the result on the existence of the smooth density for $X(t)$, which have already stated in Remark \ref{Rem:density of SDDE X(t)}, under the uniformly elliptic condition \eqref{condition:UE}. 
Since a continuous function with linear growth order can be approximated by bounded continuous functions, we can extend to the class $\mathfrak{F} (\mathbb{R} \, ; \, \mathbb{R})$ easily. 
\hfill $\square$

\begin{Rem}
{\rm{
From Proposition 1.3.5 (p.40) in \cite{Nualart}, it holds that 
\begin{align*}
& \delta \left( \frac{U(\cdot) \, \Lambda (t)}{A_1 \big( X(\cdot -r) \big) \, X( \cdot )} \, \mathbb{I}_{[0 \vee (t-r),t]} ( \cdot ) \right) \\
& = \delta \left( \frac{U( \cdot )}{A_1 \big( X(\cdot -r) \big) \, X( \cdot )} \, \mathbb{I}_{[0 \vee (t-r),t]} ( \cdot ) \right) \, \Lambda (t) 
- \int_{0 \vee (t-r)}^t \frac{U( u )}{A_1 \big( X(u-r) \big) \, X( u )} \, \frac{ \mbox{d}}{\mbox{d}u} D_u \Lambda (t) \, \mbox{d}u \\
& = \left( \int_{0 \vee (t-r)}^t \frac{U( u )}{A_1 \big( X(u-r) \big) \, X( u )} \, \mbox{d}W(u) \right) \, \Lambda (t) 
- \int_{0 \vee (t-r)}^t \frac{U( u )}{A_1 \big( X(u-r) \big) \, X( u )} \, \frac{ \mbox{d}}{\mbox{d}u} D_u \Lambda (t) \, \mbox{d}u. 
\end{align*}
In particular, consider the case $0 \le t \le r$. 
Since 
$$
U(u) = \frac{ X(u)}{x}, \quad 
\Lambda (t) = 1 + x \, A_1^\prime ( x ) \, W(t) - x \, A_1^\prime ( x ) \, A_1 ( x ) \, t
$$
for $0 \le u \le t$, we can get 
\begin{align*}
\delta \left( \frac{U( \cdot ) \, \Lambda (t)}{A_1 \big( X(\cdot -r) \big) \, X( \cdot )} \, \mathbb{I}_{[0 \vee (t-r),t]} ( \cdot ) \right) 
= \frac{W(t)}{x \, A_1 ( x )} + \frac{A_1^\prime (x )}{A_1(x)} \big( W(t)^2 -t \big) - A_1 (x) \, t \, W(t). 
\end{align*}
\hfill $\square$
}}
\end{Rem}

We shall compute the Delta, that is, a kind of the Greeks for the Asian-type option associated with an asset price model with delay in the initial point $x > 0$. 
Write $\tilde{x} = \log x$. 
For $\tilde{y} \in \mathbb{R}$, define the $\mathbb{R}$-valued processes $\tilde{X}(t) = \big\{ \tilde{X}(t) \, ; \, -r \le t \le T \big\}$ and $\tilde{Y} = \big\{ \tilde{Y}(t) \, ; \, -r \le t \le T \big\}$ by
\begin{align}
& \tilde{X}(t) = 
\begin{cases}
\tilde{x} & ( -r \le t \le 0), \\
\displaystyle \tilde{x} - \frac{ \ 1 \ }{2} \int_0^t \tilde{A}_1 \big( \tilde{X}(s-r) \big)^2 \, \mbox{d}s + \int_0^t \tilde{A}_1 \big( \tilde{X}(s-r) \big) \, \mbox{d}W(s) & (0 \le t \le T), 
\end{cases}
\label{eq: special case of SDDE tilde{X}}\\[7pt]
& \tilde{Y}(t) = 
\begin{cases}
\tilde{y} & (-r \le t \le 0), \\
\displaystyle \tilde{y} + \int_0^t \exp \big( \tilde{X}(s) \big) \, \mbox{d}s & (0 \le t \le T), 
\end{cases}
\label{eq: special case of tilde{Y}}
\end{align}
where $\tilde{A}_1 (\tilde{z}) = A_1 ( \exp(z) )$. 
Remark that $X(t) = \exp \big\{ \tilde{X}(t) \big\}$ and $Y_t = \tilde{Y}(t) / t$. 
Consider the $\mathbb{R}^2$-valued process $\mathfrak{X} = \big\{ \mathfrak{X}(t) = \big( \tilde{X}(t), \tilde{Y}(t) \big) \, ; \, -r \le t \le T \big\}$ given by  $\mathfrak{X}(t) = (\tilde{x}, \tilde{y})^\ast = : \tilde{\mathfrak{X}}$ for $-r \le t \le 0$, and 
\begin{align*}
\mathfrak{X}(t) 
& = 
\begin{pmatrix}
\tilde{x} \\
\tilde{y}
\end{pmatrix}
+ \int_0^t 
\begin{pmatrix}
- \tilde{A}_1 \big( \tilde{X}(s-r) \big)^2 /2 \\
\exp \big( \tilde{X}(s) \big)
\end{pmatrix}
\, \mbox{d}s 
+ \int_0^t 
\begin{pmatrix}
\tilde{A}_1 \big( \tilde{X}(s-r) \big) \\
0 
\end{pmatrix}
\, \mbox{d}W(s) \\[5pt]
& =: 
\tilde{\mathfrak{X}} 
+ \int_0^t \mathfrak{A}_0 \big( \mathfrak{X}(s-r), \mathfrak{X}(s)\big) \, \mbox{d}s 
+ \int_0^t \mathfrak{A}_1 \big( \mathfrak{X}(s-r) \big) \, \mbox{d}W(s)
\end{align*}
for $0 \le t \le T$. 
Let $\pi : \mathbb{R}^2 \to \mathbb{R}$ be the canonical projection defined by $\pi(\mathfrak{X}) = y$ for $\mathfrak{X} = (x, y)^\ast \in \mathbb{R}^2$. 
Remark that our main interest is to study the sensitivity of 
$$
\mathbb{E} \Big[ \Phi \big( \tilde{Y}(t) / t \big) \Big] 
= \mathbb{E} \Big[ (\Phi \circ \pi) \big( \mathfrak{X}(t) / t \big) \Big] 
$$
for a certain payoff function $\Phi$. 
Before introducing our result, we shall prepare some notations. 
For $0 \le u \le T$, let $\mathfrak{Z}( \cdot, u) = \big\{ \mathfrak{Z}(t,u) \, ; \, -r \le t \le T \big\}$ be the $\mathbb{R}^2 \otimes \mathbb{R}^2$-valued process determined by the equation 
\begin{equation}
\mathfrak{Z}(t,u) = 
\begin{cases}
\bm{0} & (-r \le t \le 0 \ \mbox{or} \ t < u), \\
\displaystyle I_2 + \int_u^t \partial_1 \mathfrak{A}_1 \big( \mathfrak{X}(s-r) \big) \, \mathfrak{Z}(s-r,u) \, \mbox{d}W(s)  \\[7pt]
\quad \displaystyle  + \int_u^t \Big\{ \partial_1 \mathfrak{A}_0 \big( \mathfrak{X}(s-r), \mathfrak{X}(s) \big) \, \mathfrak{Z} ( s-r, u) \\
\qquad \qquad \displaystyle + \partial_2 \mathfrak{A}_0 \big( \mathfrak{X}(s-r), \mathfrak{X}(s) \big) \, \mathfrak{Z}(s,u) \Big\} \, \mbox{d}s & (u \le t \le T), 
\end{cases}
\end{equation}
where $I_2 \in \mathbb{R}^2 \otimes \mathbb{R}^2$ is the identity. 

\begin{Thm}[Greeks computation of the Asian option]
For $\Phi \in \mathfrak{F} (\mathbb{R}, \, ; \, \mathbb{R})$, it holds that 
\begin{equation}
\partial_x \mathbb{E} \big[ \Phi \big( \tilde{Y}(t)/t \big) \big] 
= \mathbb{E} \big[ \Phi \big( \tilde{Y}(t) / t \big) \, \Gamma_A(t) \big], 
\end{equation}
where $\mathfrak{V}(t)$ is the Malliavin covariance matrix for $\mathfrak{X}(t)$, and 
$$
\Gamma_A (t) = \frac{ \ 1 \ }{x} \, \delta \left(  \big( \mathfrak{Z}(t,\cdot) \, \mathfrak{A}_1 \big( \mathfrak{X}(\cdot -r) \big) \big)^\ast \, \mathfrak{V}(t)^{-1} \, \partial_\mathfrak{X} \mathfrak{X}(t) \right) \, 
(1,0)^\ast. 
$$
\end{Thm}
{\it{Proof}}. 
Similarly to Theorem \ref{Thm:European Delta}, the proof is based upon the standard density argument as seen in \cite{Kawai-Takeuchi, Takeuchi}. 
We shall remark that, as for the case that $\Phi$ is the indicator function, the most crucial point is the existence of the smooth density for $Y(t)$ as stated in Theorem \ref{Thm:smooth density for Y(t)}. 

We have only to discuss the case of $\Phi \in C_b^1 (\mathbb{R} \, ; \, \mathbb{R})$. 
First of all, we shall remark that 
$$
\partial_{\tilde{x}} \mathbb{E} \Big[ \Phi \big( \tilde{Y}(t) / t \big) \Big] 
= \partial_{\mathfrak{X}} \mathbb{E} \Big[ ( \Phi \circ \pi ) \big( \mathfrak{X}(t) / t \big) \Big] \, 
(1,0)^\ast
. 
$$
We can compute the Malliavin derivative of $\mathfrak{X}(t)$ as follows: 
$$
\frac{ \mbox{d}}{\mbox{d}u} D_u \mathfrak{X}(t) 
= \mathfrak{Z}(t,u) \, \mathfrak{A}_1 \big( \mathfrak{X}(u-r) \big) \, \mathbb{I}_{[0,t]} (u), 
$$
so, the Malliavin covariance matrix $\mathfrak{V}(t) := \langle D \mathfrak{X}(t), D \mathfrak{X}(t) \rangle_{\mathbb{H}_0^1}$ for $\mathfrak{X}(t)$ is 
$$
\mathfrak{V}(t) 
= \int_0^t \mathfrak{Z}(t,u) \, \mathfrak{A}_1 \big( \mathfrak{X}(u-r) \big) \, \big( \mathfrak{Z}(t,u) \, \mathfrak{A}_1 \big( \mathfrak{X}(u-r) \big) \big)^\ast \, \mbox{d}u
$$
Suppose that the inverse of $\det \mathfrak{V}(t)$ is in $\mathbb{L}^p (\Omega)$ for any $p > 1$. 
Since 
\begin{align*}
\partial_{\tilde{\mathfrak{X}}} \big\{ (\Phi \circ \pi)( \mathfrak{X}(t)) \big\} 
& = \partial (\Phi \circ \pi) (\mathfrak{X}(t)) \, \partial_{\tilde{\mathfrak{X}}} \mathfrak{X}(t) \\
& = \int_0^t \frac{ \mbox{d}}{\mbox{d}u} D_u \big\{ (\Phi \circ \pi)(\mathfrak{X}(t)) \big\} \, \big( \mathfrak{Z}(t,u) \, \mathfrak{A}_1 \big( \mathfrak{X}(u-r) \big) \big)^\ast \, \mbox{d}u \, \mathfrak{V}(t)^{-1} \, \partial_{\tilde{\mathfrak{X}}} \mathfrak{X}(t), 
\end{align*}
the integration by parts formula leads us to see that 
\begin{align*}
\partial_{\tilde{x}} \mathbb{E} \Big[ \Phi \big( \tilde{Y}(t) / t \big) \Big] 
& = \partial_{\tilde{\mathfrak{X}}} \mathbb{E} \Big[ ( \Phi \circ \pi ) \big( \mathfrak{X}(t) / t \big) \Big] \, (1,0)^\ast \\
& = \mathbb{E} \left[ \partial (\Phi \circ \pi) \big( \mathfrak{X}(t) / t \big) \, \frac{\partial_{\tilde{\mathfrak{X}}} \mathfrak{X}(t)}{t} \right] \, (1,0)^\ast \\
& = \mathbb{E} \left[ \int_0^t \frac{ \mbox{d}}{\mbox{d}u} D_u \big\{ (\Phi \circ \pi)\big( \mathfrak{X}(t) / t \big) \big\} \, \big( \mathfrak{Z}(t,u) \, \mathfrak{A}_1 \big( \mathfrak{X}(u-r) \big) \big)^\ast \, \mbox{d}u \, \mathfrak{V}(t)^{-1} \, \partial_{\tilde{\mathfrak{X}}} \mathfrak{X}(t) \right] \, (1,0)^\ast \\
& = \mathbb{E} \left[ ( \Phi \circ \pi) \big( \mathfrak{X}(t) / t \big) \, \delta \left(  \big( \mathfrak{Z}(t,\cdot) \, \mathfrak{A}_1 \big( \mathfrak{X}(\cdot -r) \big) \big)^\ast \, \mathfrak{V}(t)^{-1} \, \partial_{\tilde{\mathfrak{X}}} \mathfrak{X}(t) \right) \right] \, (1,0)^\ast \\
& =  \mathbb{E} \left[ \Phi \big( \tilde{Y}(t) / t \big) \, \delta \left(  \big( \mathfrak{Z}(t,\cdot) \, \mathfrak{A}_1 \big( \mathfrak{X}(\cdot -r) \big) \big)^\ast \, \mathfrak{V}(t)^{-1} \, \partial_{\tilde{\mathfrak{X}}} \mathfrak{X}(t) \right) \right] \, (1,0)^\ast. 
\end{align*}
Remark that 
$$
\partial_{\mathfrak{X}} \mathfrak{X}(t) = 
\begin{cases}
\! I_2 \qquad \qquad \qquad \qquad \qquad \qquad \qquad \qquad \qquad \qquad \qquad \qquad \qquad \qquad \quad (-r \le t \le 0), \\
\! \displaystyle I_2 + \int_0^t \partial_1 \mathfrak{A}_1 \big( \mathfrak{X}(s-r) \big) \, \partial_{\mathfrak{X}} \mathfrak{X}(s-r) \, \mbox{d}W(s) \\[5pt]
\displaystyle \ + \int_0^t \Bigg\{ \!\! \partial_1 \mathfrak{A}_0 \big( \mathfrak{X}(s-r), \mathfrak{X}(s) \big) \partial_{\mathfrak{X}} \mathfrak{X}(s-r) 
+ \partial_2 \mathfrak{A}_0 \big( \mathfrak{X}(s-r), \mathfrak{X}(s) \big) \partial_{\mathfrak{X}} \mathfrak{X}(s) \!\! \Bigg\} \mbox{d}s \ \ (0 \le t \le T). 
\end{cases}
$$
Let $\mathfrak{U} = \big\{ \mathfrak{U}(t) \, ; \, -r \le t \le T \big\}$ and $\hat{\mathfrak{U}} = \big\{ \hat{\mathfrak{U}}(t) \, ; \, -r \le t \le T \big\}$ be $\mathbb{R}^2$-valued processes determined by the following ordinary differential equations: 
\begin{align*}
\mathfrak{U}(t) & = \hat{\mathfrak{U}}(t) = I_2 \qquad \qquad \qquad \qquad \qquad \qquad \ \ (-r \le t \le 0), \\
\mathfrak{U}(t) & = I_2 + \int_0^t \partial_2 \mathfrak{A}_0 \big( \mathfrak{X}(s-r), \mathfrak{X}(s) \big) \, \mathfrak{U}(s) \, \mbox{d}s \qquad (0 \le t \le T), \\
\hat{\mathfrak{U}}(t) & = I_2 - \int_0^t \hat{\mathfrak{U}}(s) \, \partial_2 \mathfrak{A}_0 \big( \mathfrak{X}(s-r), \mathfrak{X}(s) \big) \, \mbox{d}s \qquad (0 \le t \le T). 
\end{align*}
Then, it is clear that $\mathfrak{U}(t) \, \hat{\mathfrak{U}}(t) = \hat{\mathfrak{U}}(t) \, \mathfrak{U}(t) = I_2$, via the It\^o formula. 

Now, we shall consider the negative-order moment of $\det \mathfrak{V}(t)$. 
Moreover, for $0 \vee (t-r) \le u \le t$, we see that $\mathfrak{Z}(t,u) = \mathfrak{U}(t) \, \hat{\mathfrak{U}}(u)$ is invertible, because 
\begin{align*}
\mathfrak{Z}(t,u) 
& = I_2 + \int_u^t \partial_1 \mathfrak{A}_1 \big( \mathfrak{X}(s-r), \mathfrak{X}(s) \big) \, \mathfrak{Z} ( s-r,u ) \, \mbox{d}W(s)  \\
& \quad + \int_u^t \Big\{ \partial_1 \mathfrak{A}_0 \big( \mathfrak{X}(s-r), \mathfrak{X}(s) \big) \, \mathfrak{Z} ( s-r, u )+ \partial_2 \mathfrak{A}_0 \big( \mathfrak{X}(s-r), \mathfrak{X}(s) \big) \, \mathfrak{Z}(s,u) \Big\} \, \mbox{d}s \\
& = I_2 + \int_u^t  \partial_2 \mathfrak{A}_0 \big( \mathfrak{X}(s-r), \mathfrak{X}(s) \big) \, \mathfrak{Z}(s,u) \, \mbox{d}s
\end{align*}
for $0 \vee (t-r) \le u \le t$, 
$$
\mathfrak{U}(t) \, \hat{\mathfrak{U}}(u) 
= I_2 + \int_u^t \partial_2 \mathfrak{A}_0 \big( \mathfrak{X}(s-r), \mathfrak{X}(s) \big) \, \mathfrak{U}(s) \, \hat{\mathfrak{U}}(u) \, \mbox{d}s,
$$
and the uniqueness of the solutions. 
Remark that 
\begin{align*}
\mathfrak{V}(t) 
& = \int_0^t \Big( \mathfrak{Z}(t,u) \, \mathfrak{A}_1 \big( \mathfrak{X}(u-r) \big) \Big) \, \Big( \mathfrak{Z}(t,u) \, \mathfrak{A}_1 \big( \mathfrak{X}(u-r) ) \big) \Big)^\ast \, \mbox{d}u \\
& \ge \int_{0 \vee (t - r)}^t \Big( \mathfrak{Z}(t,u) \, \mathfrak{A}_1 \big( \mathfrak{X}(u-r) \big) \Big) \, \Big( \mathfrak{Z}(t,u) \, \mathfrak{A}_1 \big(\mathfrak{X}(u-r) \big) \Big)^\ast \, \mbox{d}u \\
& = \mathfrak{U}(t) \, \left\{ \int_{0 \vee (t - r)}^t \Big( \hat{\mathfrak{U}}(u) \, \mathfrak{A}_1 \big( \mathfrak{X}(u-r) \big) \Big) \, \Big(  \hat{\mathfrak{U}}(u) \, \mathfrak{A}_1 \big( \mathfrak{X}(u-r) \big) \Big)^\ast \, \mbox{d}u \right\} \, \mathfrak{U}(t)^\ast \\
& = : \mathfrak{U}(t) \, \check{\mathfrak{V}}(t) \, \mathfrak{U}(t)^\ast. 
\end{align*}
Here, the second inequality is in the matrix sense. 
It is clear that $\sup_{t \in [-r,T]} \big\| \mathfrak{U} (t) \big\|$ is in $\mathbb{L}^p (\Omega)$ for any $p > 1$. 
Hence, we have to attack the negative-order moment of $\det \check{\mathfrak{V}}(t)$. 

Remark that $\mathfrak{A}_1 \big( \mathfrak{X}(u-r) \big)$ is $\mathcal{F}_{0 \vee (t-r)}$-measurable for $0 \vee (t-r) \le u \le t$. 
When we study the lower estimate of $\check{\mathfrak{V}}(t)$, we can regard the term $\mathfrak{A}_1 \big( \mathfrak{X}(u-r) \big)$ as the constant in the integrand of $\check{\mathfrak{V}}(t)$, by taking the conditional expectation on $\mathcal{F}_{0 \vee (t-r)}$. 
Let $(\alpha_0, \beta_0) \in \mathbb{R}^2$ be fixed. 
Since 
$$
\mathfrak{A}_1 \big(\alpha_0, \beta_0 \big) = 
\begin{pmatrix}
\tilde{A}_1 ( \alpha_0 ) \\
0
\end{pmatrix}
, \quad 
\mathfrak{A}_0 \big( (\alpha_0, \beta_0), (\alpha, \beta) \big) = 
\begin{pmatrix}
- \tilde{A}_1 ( \alpha_0 )^2 /2 \\
\exp (\alpha)
\end{pmatrix}
$$
for $(\alpha, \beta) \in \mathbb{R}^2$, we see that 
\begin{align*}
\big[ \mathfrak{A}_0, \mathfrak{A}_1 \big] \big( (\alpha_0, \beta_0), (\alpha, \beta) \big) 
& = - \partial_{(\alpha, \beta)} \mathfrak{A}_0 \big( (\alpha_0, \beta_0), (\alpha, \beta) \big) \, \mathfrak{A}_1 \big( \alpha_0, \beta_0 \big) \\
& = 
\begin{pmatrix}
0 \\
- \tilde{A}_1 ( \alpha_0) \, \exp (\alpha)
\end{pmatrix}
. 
\end{align*}
The dimension of the linear space generated by $\mathfrak{A}_1 \big( \alpha_0, \beta_0 \big)$ and $\big[ \mathfrak{A}_0, \mathfrak{A}_1 \big] \big( (\alpha_0, \beta_0), (\alpha, \beta) \big)$ is $2$ for all $(\alpha, \beta) \in \mathbb{R}^2$, because of the uniformly elliptic condition \eqref{condition:UE} on $A_1$. 
Hence, we can conclude that the inverse of $\det \check{\mathfrak{V}}(t)$ is in $\mathbb{L}^p (\Omega)$ for any $p > 1$, which implies that the probability law of the $\mathbb{R}^2$-valued random variable $\mathfrak{X}(t)$ admits a smooth density with respect to the Lebesgue measure over $\mathbb{R}^2$. 
See \cite{KomatsuTakeuchi} on details. 
Therefore, we can justify the existence of the smooth density for the probability law of the $\mathbb{R}$-valued random variable $Y(t) = \tilde{Y}(t) / t$ with respect to the Lebesgue measure over $\mathbb{R}$. 
The proof is complete. 
\hfill $\square$

\begin{Rem}\label{Rem: small time}
{\rm{
Consider the case of $0 \le t \le r$. 
Then, we can derive 
\begin{align*}
& 
\mathfrak{X}(t) = 
\begin{pmatrix}
\tilde{X}(t) \\
\tilde{Y}(t)
\end{pmatrix}
= 
\begin{pmatrix}
\tilde{x} \\
\tilde{y} 
\end{pmatrix}
+ \int_0^t 
\begin{pmatrix}
- \tilde{A}_1 ( \tilde{x} )^2 /2 \\
\exp \big( \tilde{X}(s) \big)
\end{pmatrix}
\, \mbox{d}s 
+ \int_0^t 
\begin{pmatrix}
\tilde{A}_1 ( \tilde{x} ) \\
0 
\end{pmatrix}
\, \mbox{d}W(s), \\ 
& \partial_{\tilde{\mathfrak{X}}} \mathfrak{X}(t) \, (1,0)^\ast 
= 
\begin{pmatrix}
\partial_{\tilde{x}} \tilde{X}(t) \\
\partial_{\tilde{x}} \tilde{Y}(t) 
\end{pmatrix} 
= 
\begin{pmatrix}
1 - \tilde{A}_1 ( \tilde{x} ) \, \tilde{A}_1^\prime ( \tilde{x} ) \, t + \tilde{A}_1^\prime ( \tilde{x} ) \, W(t) \\
\displaystyle \int_0^t \left( 1 - \tilde{A}_1 ( \tilde{x} ) \, \tilde{A}_1^\prime ( \tilde{x} ) \, s + \tilde{A}_1^\prime ( \tilde{x} ) \, W(s) \right) \, \exp \big( \tilde{X}(s) \big) \, \mbox{d}s 
\end{pmatrix}
, \\
& \mathfrak{U}(t) = 
\begin{pmatrix}
1 & 0 \\
\tilde{Y} (t) & 1
\end{pmatrix}
, \quad 
\hat{\mathfrak{U}}(t) = 
\begin{pmatrix}
1 & 0 \\
- \tilde{Y} (t) & 1
\end{pmatrix}
, \quad 
\mathfrak{Z}(t,u) = 
\begin{pmatrix}
1 & 0 \\
\tilde{Y}(t) - \tilde{Y}(u) & 1 
\end{pmatrix}
\, \mathbb{I}_{(u \le t)}. 
\end{align*}
Since $\hat{\mathfrak{U}}(u) \, \mathfrak{A}_1 (\mathfrak{X}(u-r)) = \big( \tilde{A}_1 (\tilde{x}), - \tilde{A}_1 (\tilde{x}) \, \tilde{Y}(u) \big)^\ast$, we have 
$$
\delta \Big( \big( \hat{\mathfrak{U}}(\cdot) \, \mathfrak{A}_1 (\mathfrak{X} ( \cdot - r)) \big)^\ast \Big) 
= \int_0^t\big( \tilde{A}_1 (\tilde{x}), - \tilde{A}_1 (\tilde{x}) \, \tilde{Y}(u) \big) \, \mbox{d}W(u). 
$$
Recall $\check{\mathfrak{V}}(t) = \hat{\mathfrak{U}}(t) \, \mathfrak{V}(t) \, \hat{\mathfrak{U}}(t)^\ast$. 
Then, we see that 
\begin{align*}
& \check{\mathfrak{V}}(t) 
= 
\begin{pmatrix}
\tilde{A}_1 ( \tilde{x} )^2  \, t & \displaystyle - \tilde{A}_1 ( \tilde{x} )^2 \int_0^t \tilde{Y}(v) \, \mbox{d}v \\[7pt]
\displaystyle - \tilde{A}_1 ( \tilde{x} )^2 \int_0^t \tilde{Y}(v) \, \mbox{d}v & \displaystyle \tilde{A}_1 ( \tilde{x} )^2 \int_0^t \tilde{Y}(v)^2 \, \mbox{d}v
\end{pmatrix}
, \\
& \check{\mathfrak{V}}(t)^{-1} 
= \frac{1}{\det \check{\mathfrak{V}}(t)} \, 
\begin{pmatrix}
\displaystyle \tilde{A}_1 ( \tilde{x} )^2 \int_0^t \tilde{Y}(v)^2 \, \mbox{d}v & \displaystyle \tilde{A}_1 ( \tilde{x} )^2 \int_0^t \tilde{Y}(v) \, \mbox{d}v \\[7pt]
\displaystyle \tilde{A}_1 ( \tilde{x} )^2 \int_0^t \tilde{Y}(v) \, \mbox{d}v & \tilde{A}_1 ( \tilde{x} )^2 \, t
\end{pmatrix}
= : \frac{\Theta (t)}{\det \check{\mathfrak{V}}(t)}. 
\end{align*}
Since the Malliavin derivatives of $\tilde{X}(t)$, $\partial_{\tilde{x}} \tilde{X}(t)$, $\tilde{Y}(t)$ and $\partial_{\tilde{x}} \tilde{Y}(t)$ can be computed as 
\begin{align*}
& \frac{\mbox{d}}{\mbox{d}u} D_u \tilde{X}(t) 
= \tilde{A}_1 ( \tilde{x} ) \, \mathbb{I}_{(u \le t)}, \quad 
\frac{\mbox{d}}{\mbox{d}u} D_u \partial_{\tilde{x}} \tilde{X}(t) 
= \tilde{A}_1^\prime ( \tilde{x} ) \, \mathbb{I}_{(u \le t)}, \\
& \frac{\mbox{d}}{\mbox{d}u} D_u \tilde{Y}(t) 
= \tilde{A}_1 ( \tilde{x} ) \, \big( \tilde{Y}(t) - \tilde{Y}(u) \big) \, \mathbb{I}_{(u \le t)}, \\
& \frac{\mbox{d}}{\mbox{d}u} D_u \partial_{\tilde{x}} \tilde{Y}(t) 
= \tilde{A}_1^\prime (\tilde{x}) \, \big( \tilde{Y}(t) - \tilde{Y}(u) \big) \, \mathbb{I}_{(u \le t)} + \tilde{A}_1 (\tilde{x}) \, \big( \partial_{\tilde{x}} \tilde{Y}(t) - \partial_{\tilde{x}} \tilde{Y}(u) \big) \, \mathbb{I}_{(u \le t)}, 
\end{align*}
we see that 
\begin{align*}
& \frac{\mbox{d}}{\mbox{d}u} D_u \det \check{\mathfrak{V}}(t) 
= 2 \, \tilde{A}_1 (\tilde{x})^5 \left( \int_0^t \!\! \int_u^t \big( \tilde{Y}(v) - \tilde{Y}(u) \big) \, \big( \tilde{Y}(v) - \tilde{Y}(\sigma) \big) \, \mbox{d} v \, \mbox{d} \sigma \right) \, \mathbb{I}_{(u \le t)}, \\[5pt]
& \frac{\mbox{d}}{\mbox{d}u} D_u \Theta (t)
= \tilde{A}_1 (\tilde{x})^3
\begin{pmatrix}
\displaystyle 2 \int_u^t \tilde{Y}(v) \, \big( \tilde{Y}(v) - \tilde{Y}(u) \big) \, \mbox{d}v  & \displaystyle \int_u^t \big( \tilde{Y}(v) - \tilde{Y}(u) \big) \, \mbox{d}v \\[7pt]
\displaystyle\int_u^t \big( \tilde{Y}(v) - \tilde{Y}(u) \big) \, \mbox{d}v & 0
\end{pmatrix}
\, \mathbb{I}_{(u \le t)}, \\[5pt]
& \frac{\mbox{d}}{\mbox{d}u} D_u \hat{\mathfrak{U}}(t) 
= \tilde{A}_1 (\tilde{x}) 
\begin{pmatrix}
0 & 0 \\
- \big( \tilde{Y}(t) - \tilde{Y}(u) \big) & 0 
\end{pmatrix}
\, \mathbb{I}_{(u \le t)}, \\[5pt]
& \frac{\mbox{d}}{\mbox{d}u} D_u \partial_{\tilde{\mathfrak{X}}} \mathfrak{X}(t) \, (1,0)^\ast 
= \partial_{\tilde{x}} 
\begin{pmatrix}
\tilde{A}_1 ( \tilde{x} ) \\
\tilde{A}_1 (\tilde{x}) \, \big( \tilde{Y}(t) - \tilde{Y}(u) \big)
\end{pmatrix}
\, \mathbb{I}_{(u \le t)}. 
\end{align*}
Write $\Delta (t) := t \int_0^t \tilde{Y}(v)^2 \, \mbox{d}v - \left( \int_0^t \tilde{Y}(v) \, \mbox{d}v \right)^2$. 
Since $\mathfrak{Z}(t,u) = \mathfrak{U}(t) \, \hat{\mathfrak{U}}(u)$ for $0 \le t \le r$, Proposition 1.3.5 (p.40) in \cite{Nualart} enables us to obtain 
\begin{align*}
& \delta \Big( \big( \mathfrak{Z}(t,\cdot) \, \mathfrak{A}_1 (\mathfrak{X} ( \cdot - r)) \big)^\ast \, \mathfrak{V}(t)^{-1} \, \partial_{\tilde{\mathfrak{X}}} \mathfrak{X}(t) \Big) \, (1,0)^\ast \\
& = \delta \Big( \big( \hat{\mathfrak{U}}(\cdot) \, \mathfrak{A}_1 (\mathfrak{X} ( \cdot - r)) \big)^\ast \Big) \, \frac{\Theta (t)}{\det \check{\mathfrak{V}}(t)} \, \hat{\mathfrak{U}}(t) \, \partial_{\tilde{\mathfrak{X}}} \mathfrak{X}(t) \, (1,0)^\ast \\
& \quad - \int_0^t \Big( \hat{\mathfrak{U}}(u) \, \mathfrak{A}_1 (\mathfrak{X} ( u- r)) \Big)^\ast \, \frac{\mbox{d}}{\mbox{d}u} D_u \left( \frac{\Theta (t)}{\det \check{\mathfrak{V}}(t)} \, \hat{\mathfrak{U}}(t) \, \partial_{\tilde{\mathfrak{X}}} \mathfrak{X}(t)  \, (1,0)^\ast \right) \, \mbox{d}u \\
& = \left( \int_0^t \left( \tilde{A}_1 (\tilde{x}), - \tilde{A}_1 (\tilde{x}) \, \tilde{Y}(u) \right) \, \mbox{d}W(u) \right) \, \frac{\Theta (t)}{\det \check{\mathfrak{V}}(t)} \, \hat{\mathfrak{U}}(t) \, \partial_{\tilde{\mathfrak{X}}} \mathfrak{X}(t) \, (1,0)^\ast \\
& \quad + \int_0^t \left( \tilde{A}_1 (\tilde{x}), - \tilde{A}_1 (\tilde{x}) \, \tilde{Y}(u) \right) \, \frac{\frac{\mbox{d}}{\mbox{d}u} D_u \det \check{\mathfrak{V}}(t)}{\det \check{\mathfrak{V}}(t)^2} \, \Theta (t) \, \hat{\mathfrak{U}}(t) \, \partial_{\tilde{\mathfrak{X}}} \mathfrak{X}(t) \, (1,0)^\ast \, \mbox{d}u \\
& \quad - \int_0^t \left( \tilde{A}_1 (\tilde{x}), - \tilde{A}_1 (\tilde{x}) \, \tilde{Y}(u) \right) \, \frac{1}{\det \check{\mathfrak{V}}(t)} \, \left( \frac{\mbox{d}}{\mbox{d}u} D_u \Theta (t) \right) \, \hat{\mathfrak{U}}(t) \, \partial_{\tilde{\mathfrak{X}}} \mathfrak{X}(t) \, (1,0)^\ast \, \mbox{d}u \\
& \quad - \int_0^t \left( \tilde{A}_1 (\tilde{x}), - \tilde{A}_1 (\tilde{x}) \, \tilde{Y}(u) \right) \, \frac{\Theta (t)}{\det \check{\mathfrak{V}}(t)} \, 
\left( \frac{\mbox{d}}{\mbox{d}u} D_u \hat{\mathfrak{U}}(t) \right) \, \partial_{\tilde{\mathfrak{X}}} \mathfrak{X}(t) \, (1,0)^\ast \, \mbox{d}u \\
& \quad - \int_0^t \left( \tilde{A}_1 (\tilde{x}), - \tilde{A}_1 (\tilde{x}) \, \tilde{Y}(u) \right) \, \frac{\Theta(t)}{\det \check{\mathfrak{V}}(t)} \, \hat{\mathfrak{U}}(t) \, \left( \frac{\mbox{d}}{\mbox{d}u} D_u \partial_{\tilde{\mathfrak{X}}} \mathfrak{X}(t) \, (1,0)^\ast \right) \, \mbox{d}u \\
& = \left( \int_0^t  \left( 1, - \tilde{Y}(u) \right) \, \mbox{d}W(u) \right) \, \frac{\Theta (t)}{\tilde{A}_1 (\tilde{x})^3 \, \Delta(t)} \, \check{\mathfrak{U}}(t) \, 
\begin{pmatrix}
\partial_{\tilde{x}} \tilde{X}(t) \\
\partial_{\tilde{x}} \tilde{Y}(t)
\end{pmatrix}
\\
& \quad + \int_0^t \left( 1, - \tilde{Y}(u) \right) \, \frac{2 \, \int_0^t \!\! \int_u^t \big( \tilde{Y}(v) - \tilde{Y}(u) \big) \, \big( \tilde{Y}(v) - \tilde{Y}(\sigma) \big) \, \mbox{d} v \, \mbox{d} \sigma}{\tilde{A}_1 (\tilde{x})^2 \, \Delta(t)^2} \, \mbox{d}u \, 
\Theta (t) \, \check{\mathfrak{U}}(t) \, 
\begin{pmatrix}
\partial_{\tilde{x}} \tilde{X}(t) \\
\partial_{\tilde{x}} \tilde{Y}(t)
\end{pmatrix} 
\\
& \quad - \int_0^t \left( 1, - \tilde{Y}(u) \right) \, \frac{1}{\Delta(t)} \, 
\begin{pmatrix}
\displaystyle \int_u^t 2 \tilde{Y}(v) \, \big( \tilde{Y}(v) - \tilde{Y}(u) \big) \, \mbox{d}v  & \displaystyle \int_u^t \big( \tilde{Y}(v) - \tilde{Y}(u) \big) \, \mbox{d}v \\[7pt]
\displaystyle \int_u^t \big( \tilde{Y}(v) - \tilde{Y}(u) \big) \, \mbox{d}v & 0
\end{pmatrix}
\, \mbox{d}u 
\, \check{\mathfrak{U}}(t) \, 
\begin{pmatrix}
\partial_{\tilde{x}} \tilde{X}(t) \\
\partial_{\tilde{x}} \tilde{Y}(t)
\end{pmatrix} 
\\
& \quad - \int_0^t \left( 1, - \tilde{Y}(u) \right) \, \frac{\Theta (t)}{\tilde{A}_1 (\tilde{x})^2 \, \Delta(t)} \, 
\begin{pmatrix}
0 & 0 \\
- \big( \tilde{Y}(t) - \tilde{Y}(u) \big) & 0 
\end{pmatrix}
\, \, \mbox{d}u \, 
\begin{pmatrix}
\partial_{\tilde{x}} \tilde{X}(t) \\
\partial_{\tilde{x}} \tilde{Y}(t)
\end{pmatrix} 
\\
& \quad - \int_0^t \left( 1, - \tilde{Y}(u) \right) \, \frac{\Theta (t) }{\tilde{A}_1 (\tilde{x})^3 \, \Delta(t)} \, \check{\mathfrak{U}}(t) \, \partial_{\tilde{x}} 
\begin{pmatrix}
\tilde{A}_1 ( \tilde{x} ) \\
\tilde{A}_1 (\tilde{x}) \, \big( \tilde{Y}(t) - \tilde{Y}(u) \big)
\end{pmatrix}
\, \mbox{d}u. 
\end{align*}
\hfill $\square$
}}
\end{Rem}

\begin{Rem}
{\rm{
Consider the case $A_1(z) = \alpha_1$ and $\tilde{y} = 0$, where $\alpha_1$ is a constant. 
Since the equations \eqref{eq: special case of SDDE tilde{X}} and \eqref{eq: special case of tilde{Y}} are 
$$
\tilde{X}(t) = \tilde{x} - \frac{ \ \alpha_1^2 \ }{2} \, t + \alpha_1 \, W(t), \quad 
\tilde{Y}(t) = \int_0^t \exp \big( \tilde{X}(s) \big) \, \mbox{d}s
$$
for $t \in [0,T]$, our situation here is the Asian-type option for the classical Black-Scholes model under a risk-neutral measure $\mbox{d} \mathbb{P}$. 
As seen in Remark \ref{Rem: small time}, we can also compute 
\begin{align*}
& \delta \Big(  \big( \mathfrak{Z}(t,\cdot) \, \mathfrak{A}_1 \big)^\ast \, \mathfrak{V}(t)^{-1} \, \partial_{\tilde{\mathfrak{X}}} \mathfrak{X}(t) \Big) \, (1,0)^\ast \\
& = \frac{1}{\alpha_1 \, \Delta(t)} \, \left( W(t), \, - \int_0^t \tilde{Y}(u) \, \mbox{d}W(u) \right) \, 
\begin{pmatrix}
\displaystyle \int_0^t \tilde{Y}(v)^2 \, \mbox{d}v \\[7pt]
\displaystyle \int_0^t \tilde{Y}(v) \, \mbox{d}v 
\end{pmatrix} 
\\
& \quad + \frac{2}{\Delta(t)^2} \, \int_0^t \left( 1, - \tilde{Y}(u) \right) \, \left( \int_0^t \!\! \int_u^t \big( \tilde{Y}(v) - \tilde{Y}(u) \big) \, \big( \tilde{Y}(v) - \tilde{Y}(\sigma) \big) \, \mbox{d} v \, \mbox{d} \sigma \right) \, \mbox{d}u
\begin{pmatrix}
\displaystyle \int_0^t \tilde{Y}(v)^2 \, \mbox{d}v \\[7pt]
\displaystyle \int_0^t \tilde{Y}(v) \, \mbox{d}v 
\end{pmatrix}
\\
& \quad - \frac{1}{\Delta(t)} \, \int_0^t \left( 1, - \tilde{Y}(u) \right) \, 
\begin{pmatrix}
\displaystyle \int_u^t 2 \, \tilde{Y}(v) \, \big( \tilde{Y}(v) - \tilde{Y}(u) \big) \, \mbox{d}v \\[7pt]
\displaystyle \int_u^t \big( \tilde{Y}(v) - \tilde{Y}(u) \big) \, \mbox{d}v 
\end{pmatrix}
\, \mbox{d}u. 
\end{align*}
\hfill $\square$
}}
\end{Rem}

\begin{Rem}
{\rm{
Let $0 < r_0 \le r$ be a constant, and denote by $\Pi : C([-r,0] \, ; \, \mathbb{R}) \to C([-r, -r_0] \, ; \, \mathbb{R})$ the projection such that $\Pi (f) = \{ f(s) \, ; \, -r \le s \le -r_0 \}$ for $f \in C([-r,0] \, ; \, \mathbb{R})$. 
Similarly to the studies stated above, we can also discuss the case where the process $X = \big\{ X(t) \, ; \, -r \le t \le T \big\}$ is determined by the equation: 
\begin{equation}
X(t) = 
\begin{cases}
x & (-r \le t \le 0) \\
\displaystyle x + \int_0^t A_0 \big( \Pi( X_s) \big) \, X(s) \, \mbox{d}s 
+ \int_0^t A_1 \big( \Pi (X_s) \big) \, X(s) \, \mbox{d}W(s) & (0 < t \le T), 
\end{cases}
\end{equation}
where $A_i \in C_{0+,b}^1 \big( C([-r, -r_0] \, ; \, \mathbb{R}) \, ; \, \mathbb{R} \big) \ (i = 0, \, 1)$ with the uniformly elliptic condition on $A_1$: 
\begin{equation}
\inf_{f \in C([-r,0] \, ; \, \mathbb{R})} \big( A_1 (\Pi(f)) \big)^2 \ge C_{13}. 
\end{equation}
\hfill $\square$
}}
\end{Rem}

\section*{Acknowledgements}
This work is supported by JSPS KAKENHI Grant Number 23740083 (Grant-in-Aid for Young Scientists (B)). 
This work was motivated by a private discussion with Alexey Kulik during his stay in Japan on February in 2015. 
The author is very grateful to him for his valuable comments and suggestions. 



}}
\end{document}